\documentclass[reqno,twoside,final]{amsart}
\usepackage{hyperref}
\usepackage{amssymb,amsmath,epsf,epsfig}

\newtheorem{Conj}{Conjecture}
\newtheorem{Thm}{Theorem}
\newtheorem{Cor}{Corollary}
\newtheorem{Lem}{Lemma}

\theoremstyle{definition}
\newtheorem{Def}{Definition}
\theoremstyle{remark}
\newtheorem{Rem}{Remark}
\newtheorem{Ex}{Example}

\numberwithin{equation}{section}

\newcommand{\thmref}[1]{Theorem~\ref{#1}}
\newcommand{\secref}[1]{\paragraph\ref{#1}}
\newcommand{\lemref}[1]{Lemma~\ref{#1}}

\def\a{\alpha}
\def\b{\beta}

\def\dist{\operatorname{dist}}

\let\paragraph=\S
\let\hataccent=\^
\def\S{\varSigma}

\def\({\left(}
\def\){\right)}
\def\[{\left[}
\def\]{\right]}

\def\k{\operatorname{\Bbbk\kern1pt}}

\def\:{\colon}
\def\~#1{\widetilde{#1}}
\def\^#1{\widehat{#1}}
\def\=#1{\check{#1}}
\let\le=\leqslant
\let\ge=\geqslant

\def\<{\left<}
\def\>{\right>}
\begin{document}
\title{Convex-concave body in $\mathbb{R}P^3$ contains a line}

\author{A. Khovanskii, D. Novikov}
\date{\today}
\address{Department of Mathematics, Purdue University, West Lafayette IN}
 \email{dmitry@math.purdue.edu}
\address{Department of Mathematics, Toronto University, Toronto,
Canada}
 \email{askold@math.toronto.edu}
\thanks{Khovanskii's work is partially supported by
Canadian Grant {\rm N~0GP0156833}. Novikov's work was supported by
the Killam grant of P. Milman and by James S. McDonnell
Foundation.}

\begin{abstract}
We define a class of $L$-convex-concave subsets of $\mathbb{R}P^3$,
where $L$ is a projective line in $\mathbb{R}P^3$. These are sets
whose sections by any plane containing $L$ are convex and
concavely depend on this plane. We prove a version of Arnold
hypothesis for these sets, namely we prove that each such set
contains a line.
\end{abstract}

\maketitle

\tableofcontents

\section{Introduction}\label{intro}
Consider a connected closed hypersurface $M$ without a boundary
embedded to ${\mathbb R}P^n$. Suppose that the  second fundamental
form  of $M$ is everywhere negatively defined. It means that in
some affine coordinates in ${\mathbb R}P^n$ the hypersurface is
locally defined as $x_n=-x_1^2-...-x_{n-1}^2+\text{\rm higher
order terms}$. A well known theorem claims then that $M$ bounds a convex body in
${\mathbb R}P^n$, i. e. doesn't
intersects some hyperplane $H\subset {\mathbb R}P^n$ and bounds a
convex body in the affine space ${\mathbb R}P^n\setminus H$.

Arnold (in \cite{Arnold88}) conjectured that an analogue of this
fact holds for any hypersurface with an everywhere non-degenerate
second fundamental form. We will say that a quadratic form in
${\mathbb R}^{n-1}$ has signature $(n-k-1,k)$ if its restriction
to some $k$-dimensional linear subspace is negatively defined and
its restriction to some $n-k-1$-dimensional linear subspace is
positively defined.
\begin{Conj}[Arnold Conjecture]
Consider a domain $U\subset{\mathbb R}P^n$  bounded by  a
connected smooth hypersurface $B$. Suppose that the second
fundamental form of $B$ is  non-degenerate at any point of $B$ and
has signature $(n-k-1, k)$ (necessarily the same for all points)
with respect to the outward normal.  Then there exist a projective
subspace $L^k$ of dimension $k$ contained in $U$ and a projective
subspace $L^{n-k-1}$ of dimension $n-k-1$ not intersecting $U$.
\end{Conj}

\begin{Ex}
Domain $U=\{(x_0,...,x_n)\in {\mathbb R}P^n|
x_0^2+...+x_k^2-x_{k+1}^2-...-x_n^2\le 0\}$, bounded by a quadric,
satisfies to the conditions and conclusions of this conjecture.
\end{Ex}
\begin{Ex}
For $k=n-1$ the conditions imposed on $B$ in the conjecture
coincide with the conditions of the theorem, and
the claim of the conjecture means existence of a hypersurface not
intersecting $U$ and of an interior point of $U$. So for $k=n-1$ the
conjecture follows from the theorem above.
\end{Ex}

\subsection{Affine version of the Arnold conjecture} There is an
affine version of the Arnold conjecture: in the statement of the
conjecture  ${\mathbb R}P^n$ is changed to ${\mathbb R}^n$ and
projective subspaces to the affine one. We prove it (in \cite
{NKh3})  for surfaces asymptotically approaching to the quadratic
cone $K=\{x_1^2+...+x_k^2-x_{k+1}^2-...-x_n^2=0\}$ as
$|x_n|\to+\infty$. This condition in particular guarantees the
smoothness of the closure of these surfaces after embedding in
${\mathbb R}P^n$.

However, in the case of slightly different asymptotical behavior
the claim is wrong already for $k=1, n=3$. Consider a union
$K'\subset {\mathbb R}^3$ of moved apart  halves of $K$ (e.g.
$K'=\{(x,y,z)\quad|\quad x^2+y^2=(|z|-1)^2, |z|\ge 1\}$). We
construct (in \cite {NKh3}) an example of a domain in ${\mathbb
R}^3$ not containing lines, satisfying conditions of the affine
version of Arnold conjecture and  which boundary  asymptotically,
as $|z|\to\infty$, approaches  $K'$.

However, the closure of such domains in ${\mathbb R}P^3$ will be
non-smooth. Moreover, it cannot be made smooth by small
perturbation without creating points of degeneracy of the second
fundamental form.

\subsection{$L$-convex-concave subsets of ${\mathbb R}P^3$}

 In this paper we prove the first nontrivial case ($k=1,n=3$) of
the Arnold  conjecture in some additional assumptions. Namely, for
any projective subspace $L\subset {\mathbb R}P^n$ we define a
class of {\em $L$-convex-concave} subsets of ${\mathbb R}P^n$.
\begin{Def}\label{defLccset}
A closed set $A\subset {\mathbb R}P^n$ is $L$-convex-concave if
\begin{enumerate}
\item $A\cap L=\emptyset$,
\item for any projective subspace $N\subset {\mathbb
R}P^n$ of dimension $\dim L+1$ and containing $L$ the intersection
$A\cap N$ is convex,
\item for any projective subspace $T\subset L$ of dimension $\dim
L-1$ the complement to the  image of $\pi(A)$ under projection
$\pi:{\mathbb R}P^n\setminus T\to {\mathbb R}P^n/T$ is an open
convex set.
\end{enumerate}
\end{Def}

In general the boundary of a $L$-convex-concave subset of
${\mathbb R}P^n$ need not be smooth, so the class of
$L$-convex-concave domains is not included into the class of
domains described in the Arnold conjecture. However, any
$L$-convex-concave set after a suitable arbitrarily small perturbation
will have a smooth and non-degenerate boundary and will  satisfy
conditions of the Arnold conjecture.

The inverse inclusion is also wrong: not all domains satisfying
the conditions of Arnold conjecture are $L$-convex-concave for
some $L$. The difference is twofold. First, in the very definition
of the $L$-convex-concave domain we postulate the existence of one
of the subspaces whose existence is claimed in the Arnold
conjecture. Second, in the definition of $L$-convex-concave
domains we suppose that all its sections by subspaces containing
$L$ as a hyperplane are convex, which is a very strong assumption.

An analogue of the Arnold conjecture for  $L$-convex-concave
domains  is the following
\begin{Conj}
Any $L$-convex-concave domain $A\subset {\mathbb R}P^n$ contains a
projective subspace of  dimension equal to $n-\dim L-1$.
\end{Conj}

In this paper we prove the first nontrivial case of this
conjecture:
\begin{Thm}\label{thm:main}
Any $L^1$-convex-concave set $A\subset {\mathbb R}P^3$, $\dim
L^1=1$, contains a projective line.
\end{Thm}

\subsection{Structure of the paper.}
The proof of this theorem belongs in fact to the realm of the
convex geometry. It  heavily exploits the two fundamental theorems
of the convex geometry: Helly theorem and the Browder theorem.
Proof  is partly guided by the general ideology of the Chebyshev
best approximation. In particular, one of the key ingredients of
the proof is an analogue of the Chebyshev alternance, see
\lemref{Chebyshev property} and \thmref{nontrivial codes}.

 Further we will consider only bodies $L$-convex-concave
with respect to some fixed once and forever real projective line
$L$. So we will use the term {\em convex-concave} for the
$L$-convex-concave bodies.

Also, we will use an equivalent definition of a convex-concave
set. Namely, in \cite{NKh2} it is shown that the convex-concave
subsets of $\mathbb{RP}^3$ can be characterized in the following
way.
\begin{Def}\label{def of proj convex-concave}
A body $\mathbb{B}\in \mathbb{RP}^3$ is called {\em projective
convex-concave} with respect to a line $L$ (further called {\em
infinite line}) not intersecting  $\mathbb{B}$ if
\begin{itemize}
\item  sections of $\mathbb{B}$
by planes passing through this line (further called {\em
horizontal planes}) are all convex and
\item for any three such horizontal sections through any point of any
of them passes a line intersecting two another.
\end{itemize}
\end{Def}

\begin{Rem} One can define an affine analogue of  projective convex-concave sets. Namely,
a body  $\mathbb{B}\in \mathbb{R}^3$ is called affine
convex-concave if, first, its horizontal sections  are all convex
and, second, for any three horizontal sections through any point
of the {\em middle} one passes a line intersecting two
another.

In \cite{NKh3} we build a counterexample to an affine version of
Arnold conjecture by smoothening a suitable affine convex-concave
body. \end{Rem}

The proof is organized as follows. In \secref{HellyBrower} we show
that it is enough to prove existence of a line intersecting any
five sections of the body, see \thmref{five is enough}. This is a
standard application of the Helly theorem. From the other hand,
using Browder theorem, we prove that for any four sections we can
find a line intersecting all of them, see \thmref{Browder}.

Starting from \secref{Chebyshev} we are dealing with five fixed
sections of a convex-concave body.  The general idea is simple.
Fix an Euclidean metric on some affine cart in ${\mathbb R}P^3$
containing all five sections and take a line closest to these five
sections (the Chebyshev line). Our goal is to prove that one  can
always find a line which lies closer to these five sections,
unless the Chebyshev line intersects all five sections.

More exact, in \secref{Chebyshev} we introduce the Euclidean
metric, define the Chebyshev line and prove its basic properties.
On planes containing sections arise five half-planes with the
property that any line  lying closer to five sections than the
Chebyshev line  should intersect all these half-planes. The
opposite is almost true. Namely, any line intersecting these
half-planes (further called {\em good deformation}) produce a line
closer to the sections than the Chebyshev line, see \lemref{good
deformation}. So all we need to prove is the existence of a line
intersecting these five half-planes, which depends on the
projective properties of their mutual position only. These
properties are the main object of further investigations.

At this stage a split occurs. We impose a condition of genericity
on the collection of these half-planes (namely, their boundaries
should be pairwise non-parallel) and deal further with
non-degenerate cases only. In degenerate cases existence of the
good deformation follows from \thmref{Browder}  due to a
remarkable self-duality of the condition of $L$-convex-concavity,
see \secref{ssec:degeneration} and \cite{NKh2}.

In \secref{codes} and \secref{boards} we investigate combinatorial
properties of a collection of five half-planes corresponding to a
Chebyshev line, forgetting for a moment the convex-concavity
condition. In other words, we consider  a more general problem of
properties of a line closest to five convex figures on five
parallel planes. This reduces to a purely combinatorial problem
about possible arrangements of rooks on a chess board. We find an
equivalent of the classical condition of Chebyshev alternance for
our situation. Namely, only six possible combinatorial types of
collections of half-planes are possible, see \thmref{nontrivial
codes}.

In \secref{final} for each of these six types we prove  existence
of a good deformation using the convex-concavity condition. More
exact, each of these combinatorial types have some continuous
parameters (e.g. distances between sections). If a configuration
of half-planes arose from a Chebyshev line, then these parameters
should satisfy some inequalities. In other words, only part of the
space of parameters corresponds to  Chebyshev alternances. It
turns out that configurations of half-planes arising from sections
of a convex-concave body belong to the complement to this part.

Namely, using the combinatorial properties of each case, we are
able to prove existence of a line intersecting four of the
half-planes in a some particular  sectors. These sectors are
chosen in such a way that the line intersecting them should
necessarily intersect the fifth half-plane and the existence of a
good deformation follows.

\section{Applications of the Helly theorem and of the Browder theorem
}\label{HellyBrower}

In this section we first introduce a linear structure on the set
of all lines not intersecting the line $L$. We prove that the
\thmref{thm:main} follows from the fact that for any five sections
of a convex-concave body there is a line intersecting all of them.
Another result claims that for any {\em four} sections there is a
line intersecting all of them.

\subsection{Linear structure on the set of all non-horizontal lines}
We will call a line {\em non-horizontal} if it doesn't intersect
the infinite line. We choose coordinates in a  complement to some
horizontal plane in such a way that the infinite line lies in the
projective plane $\{z=0\}$. In these coordinates  non-horizontal
lines have a parametrization of the type $x=az+b, y=cz+d$. This
correspondence $\{$non-horizontal line$\} \to (a,b,c,d)$ defines
coordinates on the set $\mathbb{U}$ of all non-horizontal lines.

\begin{Rem}
These coordinates are correlated with the affine structure in
horizontal planes: intersection of a convex combination of two
lines with a horizontal plane is a convex combination (with the
same coefficients) of intersections of these two lines with this
plane. Therefore the affine structure defined by these coordinates
is independent of the choice of coordinates and depends on the
choice of the infinite line only (however, the linear structure,
i.e. the line with coordinates  $(0,0,0,0)$ (=$z$-axis), can be
chosen arbitrarily).
\end{Rem}
\begin{figure}[h]
\centerline{\epsfysize=0.2\vsize\epsffile{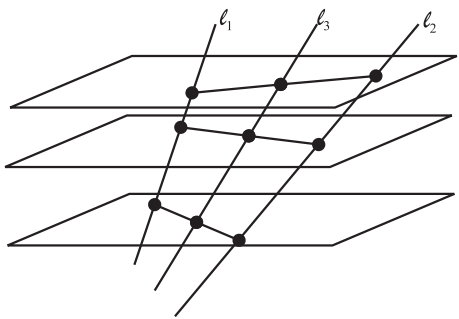}}
\caption{$\ell_3$ is a linear combination of $\ell_1$ and
$\ell_2$}
\end{figure}

Denote by $U_t$ the set of all non-horizontal lines intersecting a
horizontal section $S_t=\mathbb{B}\cap \{z=t\}$ of a projective
convex-concave body $\mathbb{B}$. From the last remark we
immediately see that
\begin{Lem}\label{conv of set of lines intersect sections}
$U_t$ is closed and is convex in the coordinates introduced above.
\end{Lem}
The inverse is also true. Namely, for any horizontal plane
$\{z=t\}$ there is a map $\phi_t:{\mathbb U}\to \{z=t\}$ mapping a
non-horizontal line to its point of intersection with this plane.
\begin{Lem}
This map preserves convexity, i.e  the image of a convex set is
again a convex set.
\end{Lem}

\subsection{Non-horizontal lines and sections of a convex-concave body}
\subsubsection{Five sections: Helly theorem}
\begin{Thm}\label{five is enough}
The \thmref{thm:main} follows from the following claim:
$$
\forall t_1,t_2,t_3,t_4,t_5\in\mathbb{R}
\qquad\bigcap_{i=1}^5U_{t_i}\not=\emptyset
$$
In other words, it is enough to prove that for any five horizontal
sections $S_i$ of $\mathbb{B}$ there exists a line intersecting
all of them.
\end{Thm}

\begin{proof}
Indeed, the \thmref{thm:main} is equivalent to
$\bigcap_tU_t\not=\emptyset$. Since $U_t$ are convex subsets of
$U\cong\mathbb{R}^4$, the claim is almost a  particular case
($n=4$) of the classical Helly theorem:

\begin{Thm}[Helly theorem, see \cite{helly,
helly2}]\label{Helly thm} Intersection of a finite family of
closed convex sets  in ${\mathbb R}^n$ is nonempty if and only if
intersection of any $n+1$ of them is nonempty.
\end{Thm}

The only problem is that the family $U_t$ is not finite. However,
one can circumvent this technicality using the fact that
\begin{Lem}
Intersection of any two different $U_t$ is compact.
\end{Lem}
Indeed, any line belonging to $U_{t_1}\cap U_{t_2}$ is uniquely
defined by its points of intersection with these two sections, so
$U_{t_1}\cap U_{t_2}$ is homeomorphic to $S_{t_1}\times S_{t_2}$,
which is compact.

So, take a compact $K=U_1\cap U_0$ and consider a family of sets
$\widetilde{U}_{t}=K\setminus U_t$. These sets are relatively open
in $K$. We  want to prove that $\cap_tU_t\not=\emptyset$. If not,
then $\widetilde{U}_{t}$ is a covering of $K$, so we can take a
finite family of $\widetilde{U}_{t_i}$ covering $K$. It means that
the intersection of a {\em finite} family consisting of the
corresponding $U_{t_i}$ and $U_1$ and $U_0$ will be empty. This is
impossible by Helly theorem if intersection of any five of $U_t$
is nonempty.
\end{proof}

\subsection{Four sections: Browder theorem}

It turns out that the convex-concavity condition (even the affine
one) guarantees existence of a line passing through any four
sections. We will prove this in slightly more general assumptions.
\begin{Thm}\label{Browder}
Let $A,B,C,D$ be four  compact convex non-empty sets in
$\mathbb{R}^n$ satisfying the following condition:
\begin{enumerate}
\item $A\subset\{x_n=t_1\}, B\subset \{x_n=t_2\}, C\subset
\{x_n=t_3\}, D\subset
\{x_n=t_4\}$, where $t_i$ are pairwise different;
\item through any point of $B$ passes a line intersecting both
$A$ and $C$, and
\item through any point of $C$ passes a line intersecting both
$B$
and $D$.
\end{enumerate}

Then there exists a line intersecting all four bodies.
\end{Thm}
\begin{Rem} Here we use only part of conditions provided by
convex-concavity.\end{Rem}

We will use a  Browder theorem --- a fixed-point theorem for upper
semi-continuous set-valued mappings, see \cite{Browder}.

Let $f:X\to Set(X)$ be a mapping from $X$ to the set of all
subsets of $X$.
\begin{Def}\label{upper-cont}
$f$ is called {\em upper semi-continuous} on $X$ if for any
$x_0\in X$ and any open set $G$ containing $f(x_0)$ there exists a
neighborhood $U$ of $x_0$ such that $f(x)\subset G$ for all $x\in
U$.
\end{Def}

\begin{Rem}
For single-valued maps this property means continuity.
\end{Rem}

Our theorem follows from  the following result of Browder:
\begin{Thm}[see \cite{Browder}]\label{Fan}
Let $X$ be a non-empty compact convex set in a real, locally
convex, Hausdorff topological vector space $E$. Let $f$ be an
upper-semicontinuous set-valued mapping defined on $X$ such that
for each $x\in X$, $f(x)$ is a non-empty closed convex subset in
$X$. Then there exists a point $\hat x\in X$ with $\hat x\in
f(\hat x)$.
\end{Thm}

We will apply this theorem to the  composition $f:B\to CSet(B)$ of
the tautological map $B\to CSet(B)$ and two maps $h_1:CSet(B) \to
CSet(C)$ and $h_2:CSet(C)\to CSet(B)$, where $CSet(B)$ and
$CSet(C)$ are sets of all compact convex subsets of $B$ and $C$
correspondingly. Namely, for $U\subset B$ we define $h_1(U)\subset
C$ as set of all points of $C$ which lie on a line intersecting
both  $A$ and $U$. Similarly, for $V\subset C$ we define
$h_2(V)\subset B$ as set of all points of $B$ which lie on a line
intersecting both $D$ and $V$. These maps are completely defined
by their restrictions to the one-point subsets of $B$ and $C$
correspondingly, namely $h_i(U)=\cup_{x\in U} h_i(\{x\})$.

Check first that our result indeed follows from the \thmref{Fan}.
Suppose that  $x\in f(x)$. It means that $x\in h_2(y)$ for some
point  $y\in h_1(\{x\})$. By definition of $h_i$ it means that the
line passing through $x$ and $y$ intersects both $A$ and $D$,
q.e.d.

We have to check that $f(x)$ satisfies conditions of
\thmref{Browder}.

By convex-concavity $f(x)$ is non-empty for all $x\in B$.

\begin{Lem}\label{f is upper semi-continuous}
$f(x)=h_2(h_1(\{x\}))$ is upper semi-continuous.
\end{Lem}

We will prove that both $h_1$ and $h_2$ are upper semi-continuous
in the sense defined below, and the claim will follow from the
fact that the composition of upper semi-continuous maps is again
upper semi-continuous. Denote by $N_{\delta}(U)=\{x| \dist(x,
U)<\delta\}$ the $\delta$-neighborhood of $U$.
\begin{Lem}\label{up-cont}
Mapping $h_1$ is upper semi-continuous in the following sense: for
any $U\in CSet(B)$ and any $\epsilon>0$ there exist a $\delta>0$
such that if $U'\subset N_{\delta}(U)$ then $h_1(U')\subset
G=N_{2\epsilon}(h_1(U))$. The mapping $h_2$ is also upper
semi-continuous.
\end{Lem}

\begin{proof}
The proof is the same for both $h_1$ and $h_2$, so we prove it for
$h_1$ only. By definition $h_1(U)=\cup_{x\in U}h_1(\{x\})$.
Therefore by compactness of $U$ it is enough to prove that for any
$b\in B$ and any $\epsilon>0$ there is a $\delta>0$ such that if
$\dist(b',b)<\delta$ then $h_1(b')\subset N_{\epsilon}(h_1(b))$.

Note that $h_1(b)=C\cap A_b$, where $A_b$ is a compact
continuously depending on $b$ in Hausdorff metric ($A_b$ and
$A_{b'}$ differ by a shift).
\begin{figure}
\centerline{\hfill\epsfysize=0.14\vsize\epsffile{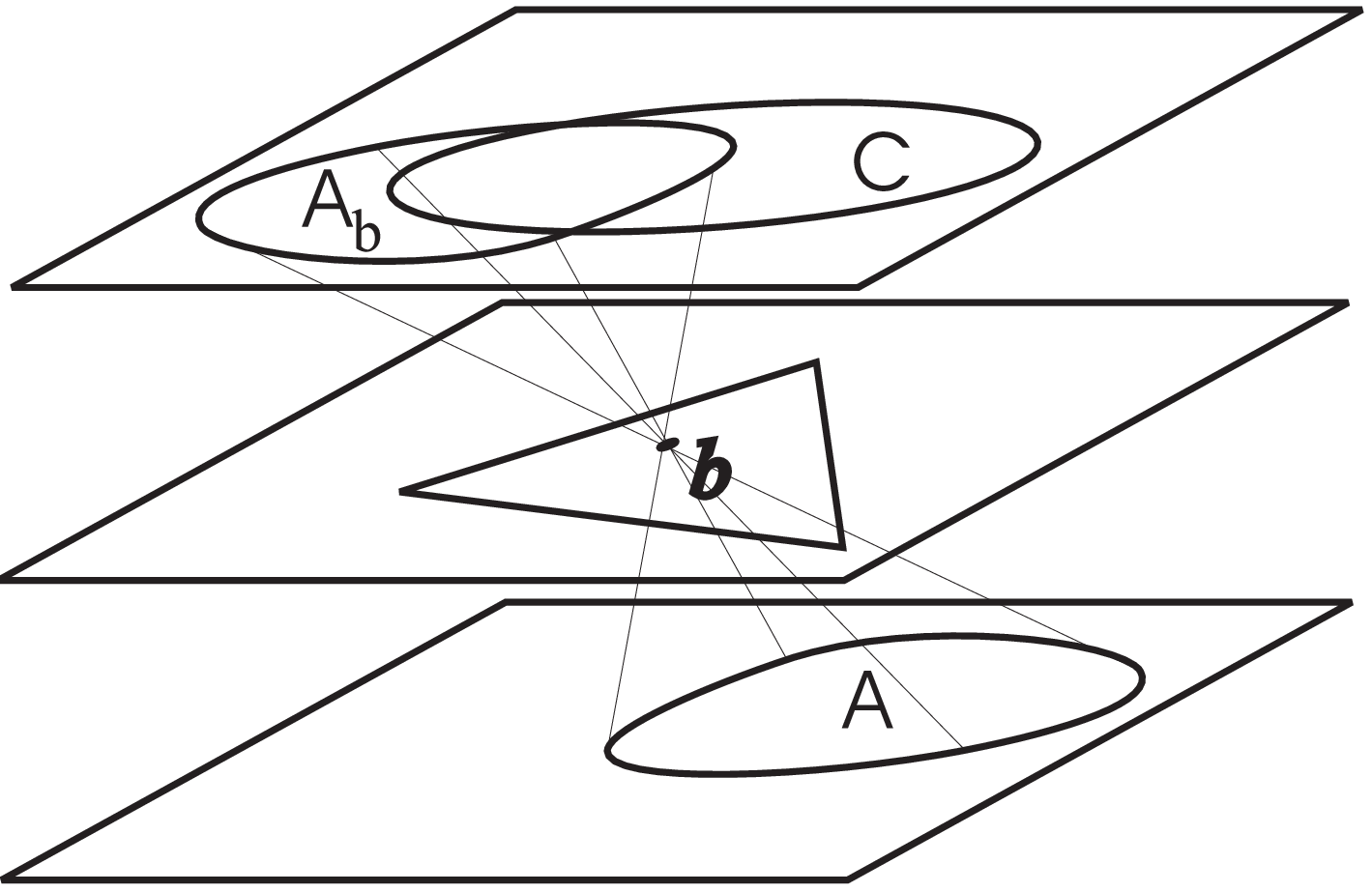}
\hfill\epsfysize=0.12\vsize\epsffile{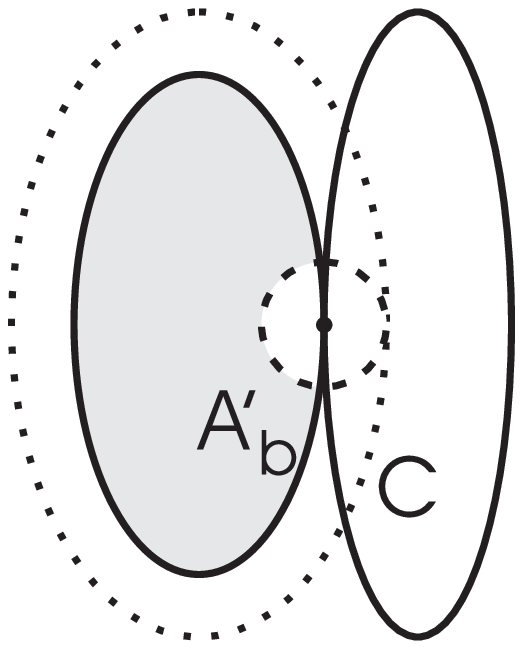}\hfill}
\caption{Construction of $h_1(b)$ and upper semi-continuity of
$h_1(U)$.}
\end{figure}

The claim follows form the fact that an intersection of a compact
with another compact continuously depending on parameters depends
upper semi-continuously on parameters.
 Let's prove this fact. Let $V=h_1(b)$,
and $A'_b=A_b\setminus N_{\epsilon}(V)$. Let
$0<\alpha<\min(\epsilon,\dist(A'_b, C))$.  For $b'$ close enough
to $b$ we have $A_{b'}\subset N_{\alpha}(A_b)$ and
$$
h_1(b')=C\cap A_{b'}\subseteq C\cap
N_{\alpha}(A_b)\subseteq(C\cap
N_{\alpha}(A'_b))\cup(C\cap
N_{\alpha}(N_{\epsilon}(V)))\subseteq
C\cap N_{\alpha+\epsilon}(V)\subset G.$$

The second inclusion is true by continuous dependence of $A_b$ on
$b$, the third is true since $A_b\subset A'_b\cup
N_{\epsilon}(V)$, the fourth is true since $C\cap
N_{\alpha}(A'_b)=\emptyset$ by choice of $\alpha$ and the last one
is true since $\alpha+\epsilon<2\epsilon$.
\end{proof}

To satisfy the last condition of the \thmref{Fan} we have to check
that $f(x)$ is a closed convex subset of $B$.
\begin{Lem}
$h_i(U)$ is compact convex set as soon as $U$ is compact
convex set.
\end{Lem}
\begin{proof}
Indeed, the set of lines intersecting both $U$ and $A$ is convex
(as intersection of two convex closed sets) and compact (since a
line is defined by its two points of intersection with $U$ and
$A$, which are both compact), so the set of points of
intersections of these lines with $\{x_n=t_3\}$ is also convex and
compact. But $h_1(U)$ is exactly the intersection of this set with
$C$, so it is also convex and compact.
\end{proof}

\begin{Rem} From a Leray theorem and the previous result we
get that the set of non-horizontal lines intersecting at least one
of the chosen five sections is homotopically equivalent to a ball
or to a sphere according to the existence or nonexistence of a
line passing through all five sections.  We know that there exist
affine convex-concave bodies (see introduction and \cite{NKh3})
without a line inside, so the case of a sphere is possible. This
sphere divides the set of all non-horizontal lines into two
connected parts. As a corollary we see that for some five sections
of these affine convex-concave body (in our example in \cite{NKh3}
these are just line segments) there is a line not intersecting
them which cannot be moved to infinity without intersecting the
sections.
\end{Rem}

\section{Chebyshev line}\label{Chebyshev}
By the previous section all we need to prove is that through any
five horizontal sections of the convex-concave body passes a line.
We fix them from now on. We choose a  sixth horizontal plane $L$
(not containing sections), choose affine coordinates in $\setminus
L\mathbb{R}^3\cong\mathbb{R}P^{3}$ and, using a standard scalar
product, introduce a metric on horizontal planes. Using this
metric we  define a Chebyshev line --- a line minimizing the
maximal distance from its point of intersection with a plane of
the section to the section. On each plane containing a section we
choose a half-plane containing the section with boundary passing
through the point of intersection of the Chebyshev line with the
plane and perpendicular to the shortest segment joining this point
to the section.

In this and the next section we investigate combinatorial
conditions imposed on the configuration of these half-planes by
the fact that the Chebyshev lines minimizes the maximal distance
to the sections.

\subsection{The Chebyshev line}

Denote by $S_1, S_2, S_3, S_4$ and $S_5$ the five sections of a
convex-concave body $B\in\mathbb{RP}^3$ cut by five horizontal
planes $L_i$, i.e. $S_i=B\cap L_i$. Choose coordinates $(\tilde
x,\tilde y,\tilde z,\tilde w)$ in $\mathbb{RP}^3$ in such a way
that the infinite line has equation $\tilde z=\tilde w=0$ and
$S_i\subset\{\tilde w\not=0\}\cong \mathbb{R}^3$. We take standard
coordinates $(x=\frac{\tilde x}{\tilde w}, y=\frac{\tilde
y}{\tilde w}, z=\frac{\tilde z}{\tilde w})$ in $\{\tilde w\not
=0\}\simeq\mathbb{R}^3$. In these coordinates the planes $L_i$ are
given by equations $L_i=\{z=t_i\}$. We take metric on $L_i$
induced by a scalar product
$$
((x_1,y_1,z_1),(x_2,y_2,z_2))=x_1x_2+y_1y_2+z_1z_2.
$$
Suppose that there is no line intersecting all five sections $S_i$
(otherwise there is nothing to prove).
\begin{Def} The  (non-horizontal) line $\ell$ minimizing the
$\max_{i=1,...,5}\dist(\ell\cap L_i, S_i)$  (where $L_i$ are the
horizontal planes containing $S_i$) will be called a Chebyshev
line.
\end{Def}
The existence of this line follows  from compactness of sections.
Further we will denote $a_i=\ell\cap L_i$ and by $s_i\in S_i$ the
point of $S_i$ closest to $a_i$.

\begin{Lem}[Chebyshev property]\label{Chebyshev property}
The $\dist(a_i, S_i)=\dist (a_i,s_i)$ are all equal.
\end{Lem}

\begin{proof}
Indeed, let one of them, say $\dist(a_1,S_1)$ is strictly smaller
than all others. By    the Browder theorem \thmref{Browder} there
exists a line $\ell_1$ which intersects all four remaining
sections. Therefore for small values of $\epsilon$ the points of
intersections of the line
$\ell_{\epsilon}=(1-\epsilon)\ell+\epsilon\ell_1$ lies closer to
$S_i$ than $a_i$ for $i=2,3,4,5$. From the other hand,
$\dist(\ell_{\epsilon}\cap L_1, S_1)$ changes continuously with
$\epsilon$. So for small $\epsilon>0$ we get
$\max_{i=1,...,5}\dist(\ell_{\epsilon}\cap L_i,
S_i)<\max_{i=1,...,5}\dist(a_i, S_i)$, which  contradicts to the
Chebyshev property of $\ell$.
\end{proof}

\begin{Cor}
The Chebyshev line $\ell$ doesn't intersect $S_i$ if there is no
line intersecting all $S_i$.
\end{Cor}

Further, in order to simplify the notations, we will suppose that
the coordinates are chosen in such a way that the Chebyshev line
coincides with the $z$ axis. Indeed, a linear transformation of
the type $(x,y,z)\to (x-(az+b), y-(cz+d),z)$ doesn't change metric
in horizontal planes, so the Chebyshev line for the shifted
sections will be the shifted Chebyshev line. From the other side,
using a transformation of this type we can move any non-horizontal
line to the $z$-axis.

\subsection{Five half-planes}
The Chebyshev condition on the line $\ell$ says that one cannot
find five points $a'_i\in L_i$ lying on a line and such that
$\dist(a'_i, S_i)< \dist(\ell\cap L_i, S_i)$. Here we describe
explicitly what the  second requirement means.

For each $a_i=\ell\cap L_i$ we can indicate an angle of desirable
directions in $L_i$: if $a_i$ moves in this direction then the
$\dist(a_i,S_i)$ decreases. These are directions forming an acute
angle with the direction $\overrightarrow{a_is_i}$. So arises the
half-plane $H_i=\{x\in L_i|(\overrightarrow{a_i
x},\overrightarrow{a_is_i})\ge 0\}$. The vector
$\overrightarrow{a_is_i}$ is orthogonal to its boundary and is
directed inward.

Another description of $H_i$ is as follows: the function
$f(x)=\dist(x,S_i)$ is a smooth function everywhere on
$L_i\setminus S_i$, so in particular  for $x=a_i$. After
identification of $T_{a_i}L_i$ and $L_i$ the half-plane $H_i$ is
described as  $\{df_{a_i}(\cdot)\ge 0\}$.

We will need further the following evident statement, see Figure
3:
\begin{Lem}\label{lem:euclid1}
Let $H$ be a half-plane in $  L_1$ bounded by a line passing
through $a_1$ and normal to $\overrightarrow{a_1n}$. Suppose that
$S_1\subset H$. Then $n\in H_1$.
\end{Lem}

\begin{figure}[h]\label{Hiandai}
\input  epsf
\centerline{\epsfysize=0.1\vsize\epsffile{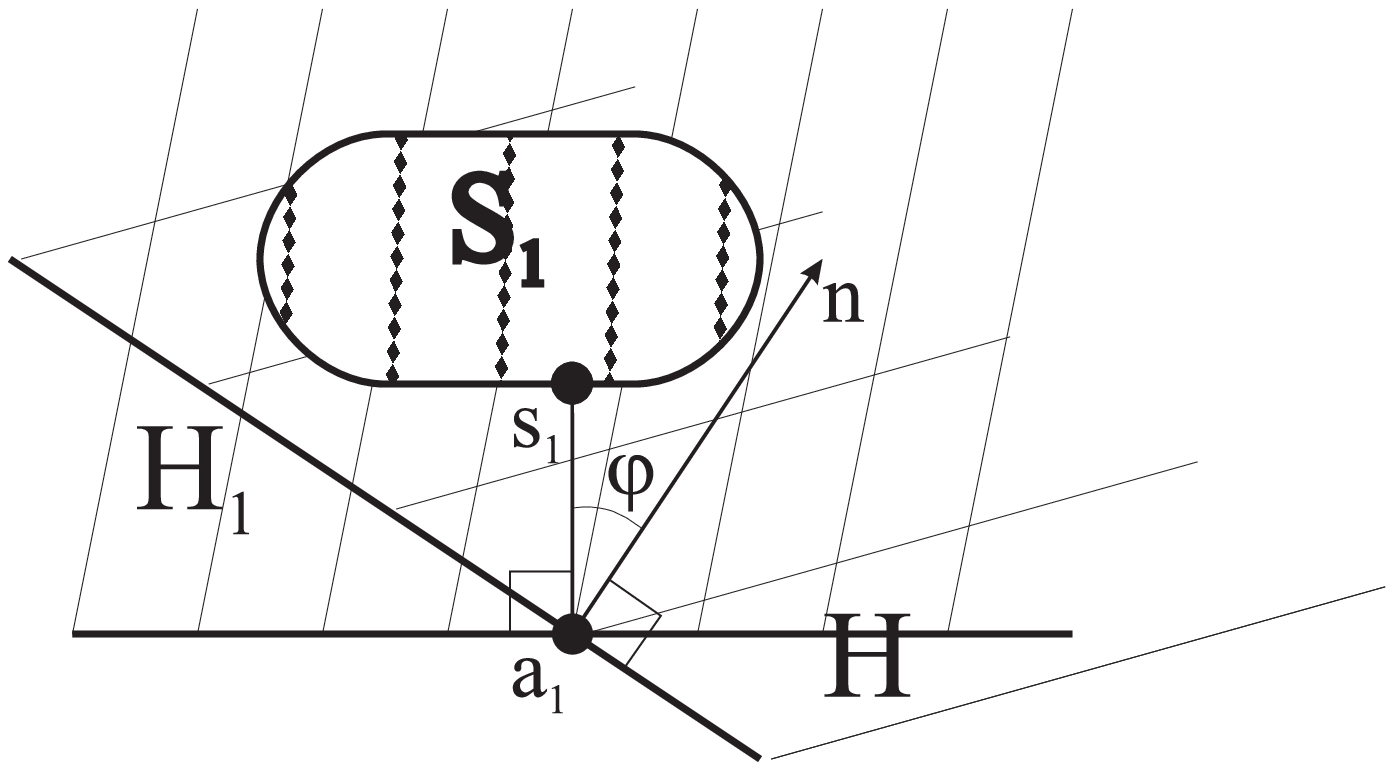}}
\caption{The half-plane $H$ contains $s_1$, so $n\in H_1$.}
\end{figure}

\subsection{Good deformations}
Here we describe lines (further called {good deformations}) whose
existence  contradicts to the fact that the Chebyshev line $\ell$
doesn't intersect the sections $S_i$ . Our goal from now is to
prove their existence.

\begin{Lem}
If $\ell$  is the Chebychev line for $S_i$ and $H_i$ are as above,
then there exist no line intersecting interiors of all $H_i$.
\end{Lem}
\begin{proof}
Suppose there exists a line $\tilde\ell$ intersecting interiors of
all $H_i$. Then $z$-axis cannot be a Chebyshev line since
$$\max_i\dist((\epsilon\tilde\ell+(1-\epsilon)\ell)\cap L_i,
S_i)<\max_i\dist(a_i, S_i)$$ for $\epsilon>0$ - small enough. In
other words, moving the Chebyshev line in the direction of
$\tilde\ell$ in the space of all non-horizontal lines decreases
its distance to $S_i$.

Indeed, all we have to check is that
$\frac{d}{d\epsilon}(\dist((\epsilon\tilde\ell+(1-\epsilon)\ell)\cap
L_i, S_i))|_{\epsilon=0}<0$, which follows directly from
definitions (of $H_i$ and the linear structure on the set of
non-horizontal lines).
\end{proof}

In fact one can prove a stronger claim.
\begin{Def}[Good deformation]\label{good deformation}
A line $\ell_1$ intersecting all $H_i$ and interior of at least
one of them will be called {\em good deformation}.
\end{Def}

\begin{Lem}\label{good deformations}
If $\ell$  is the Chebychev line for $S_i$ and $H_i$ are as above,
then there exists no  good deformation.
\end{Lem}
\begin{proof}
The proof uses the same idea as  \lemref{Chebyshev property}.
Suppose that $\ell_1$ intersects the interior of $H_1$ and denote
by $\ell_{2345}$ the line intersecting $S_2,S_3,S_4,S_5$ (it
exists by \thmref{Browder}).
 Consider the two-parametric family of
lines $\lambda\ell_1+\mu\ell_{2345}+(1-\lambda-\mu)\ell$. The idea
is that, in linear approximations, moving $\ell$ toward $\ell_1$
decreases distance to $S_1$ (while not increasing other
distances), and moving $\ell$ toward $\ell_{2345}$ decreases
distances to all other sections. So some combination of these two
movements decreases the maximal distance from the Chebyshev line
to sections, which is impossible.

In other words,  denote  points of intersection of
$\lambda\ell_1+\mu\ell_{2345}+(1-\lambda-\mu)\ell$ with $L_i$ by
$b_i^{\lambda,\mu}$. Then
$\frac{\partial}{\partial\lambda}|_{\lambda=\mu=0}\dist(b_i^{\lambda,\mu},
S_i)$ are non-negative for $i=2.3.4.5$ and is strictly positive
for $i=1$. Also
$\frac{\partial}{\partial\mu}|_{\lambda=\mu=0}\dist(b_i^{\lambda,\mu},
S_i)$ are strictly positive       for $i=2,3,4,5$. Therefore for
some positive $c_1,c_2$ we have $\frac{\partial}{\partial
\epsilon}|_{\epsilon=0}\dist(b_i^{\epsilon c_1,\epsilon c_2},
S_i)>0$ for $i=1,2,3,4,5$, i.e. for small $\epsilon>0$ the line
$\epsilon c_1\ell_1+\epsilon c_2\ell_{2345}+(1-\epsilon
c_1-\epsilon c_2)\ell$ is closer to $S_i$ than the Chebyshev line
- contradiction.
\end{proof}
\begin{Rem} The use of convex-concave property of the sections
is almost unnecessary: any four parallel half-planes with pairwise
non-parallel (see below) sides can be intersected by a line, which
is as good as $\ell_{2345}$ for the proof.
\end{Rem}

\subsection{Degenerate cases}\label{ssec:degeneration}

In what follows we will always impose the following genericity
assumption on $H_i$: we  assume that $\partial H_i$ are pairwise
non-parallel (i.e. do not intersect in $ \mathbb{R}P^3$).

 For the degenerate cases (with some of the boundaries $\partial
H_i$ being
parallel) the proof of existence of a good deformation is
reduced
via duality considerations to the \thmref{Browder}, see
\cite{NKh2}. This is done in the following way:
\begin{enumerate}
\item First, we circumscribe convex polygons $P_i$ with
$\le8$ sides around $S_i$. The sides are tangent to $S_i$ and
parallel to the boundaries of $ H_i$.
\item Second, we build the maximal (by inclusion) convex-concave
body $P$ with sections $P_i$. It exists since $S_i$ were sections
of a convex-concave body. P is the union of all points $a\in
\mathbb{R}P^3$ with the property that through any two $P_i$ and
the point $a$ passes a line.
\item Third, we consider a dual $\widetilde{P}$ of $P$ with
respect to a special duality constructed in \cite{NKh2}.
$\widetilde{P}$ is also a convex-concave body. Sections of $
\widetilde{P}$ correspond to projections of $P$. We prove that  $
\widetilde{P}$ is constructed from {\em four} convex figures in
the way described in (2). By \thmref{Browder} there exists a line
intersecting all four of them and therefore this line lies inside
$ \widetilde{P}$.
\item The dual  of this line lies inside $P$ and therefore
intersects all $P_i$. Since $P_i\subset H_i$, this line is a good
deformation.
\end{enumerate}

\section{Combinatorial properties  of half-planes arising from a Chebyshev line}\label{codes}

In this and the next chapters we investigate combinatorial
properties of mutual position of the five half-planes constructed
above. We do not use in this chapter the convex-concavity of the
sections $S_i$ (so the results are valid for any five convex
compact figures lying on five horizontal planes), and use only
part of conditions implied by the fact that $\ell$ is the
Chebyshev line for $S_i$. Namely, we use, first, the absence of
lines interior of all $H_i$ and intersecting $\ell$, and, second,
the genericity assumption of \secref{ssec:degeneration}. We single
out six combinatorial types of configurations of half-planes
satisfying these two assumptions.

The settings  we deal with can be described in projective terms.
Namely, in $\mathbb{R}P^3$ we are given a configurations
consisting of
\begin{enumerate}
\item five different projective planes $L_i$, all containing the same line (further called
infinite line),
\item five half-planes $H_i\subset L_i$ - parts of these planes - containing
convex (with respect to the infinite line) figures $S_i$ together
satisfying convex-concavity condition. Boundary of each half-plane
consists of the infinite line and some other line. The other lines
are pairwise nonintersecting by genericity assumption of
\secref{ssec:degeneration};
\item a line $\ell$ intersecting all these other lines and not intersecting the infinite line.
\end{enumerate}

In this chapter we encode  combinatorial properties of
configurations by a purely combinatorial code, leaving temporarily
aside continuous parameters of the problem (like distances between
$L_i$). This encoding can be done in several ways, so to each
configuration correspond several codes. The configurations we need
have the property that none of the corresponding codes is trivial.
In the next chapter we will see that there are at most six such
configurations.

\subsubsection{Coding}
 We will code  combinatorial properties of configurations using
projections from points $x\in\ell$ to horizontal planes. As a
result we will get a code --- a permutation of numbers $1,2,3,4,5$
with signs.

The line $\ell$ is an affine part of a projective line
$\widetilde{\ell}\cong\mathbb{RP}^1\cong\mathbb{S}^1$.  This
projective line is divided into $5$ intervals by its points of
intersection with half-planes $H_i$. We choose a point $M\in\ell$
from one of these intervals and orientation on $\ell$. We
enumerate  the points of intersections of the half-planes with
$\ell$ starting from $M$ according to the chosen orientation, thus
enumerating the half-planes by numbers $1,2,3,4,5$.

Consider a projection $\pi:\mathbb{R}^3\setminus L_M \to L_1$
(where $L_M$ is the horizontal plane passing through $M$). Take an
orientation on the circle $\mathbb{S}^1\subset L_1$ centered at
$a_1=\ell\cap L_1$ and a point $N\in \mathbb{S}^1\setminus
\cup\pi(\partial H_i)$. Thus we get an enumeration of the set of
 $10$ points  $\mathbb{S}^1\cap (\cup\pi(\partial H_i))$ (note that by
non-degeneracy assumption none of $\partial H_i$ are parallel).

We can now write down a sequence of five numbers with signs
(further called a code) which will encode the combinatorial
properties of the configuration: on the $i$-th place of this
sequence stands the number of the half-plane which boundary
projects onto the $i$-th point on $\mathbb{S}^1$ taken with $+$ if
the projection contains the point $N$ and with $-$ otherwise.

\begin{figure}[h]
\centerline{\epsfysize=2in \epsffile{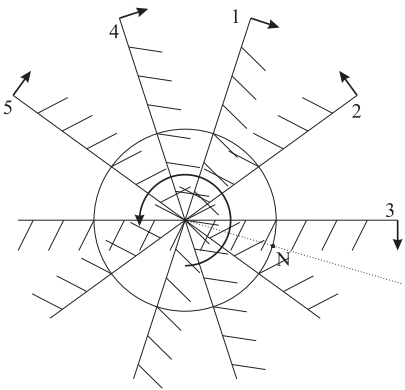}}
\caption{This projection and choices of $N$ and of orientation of
$\mathbb{S}^1$ correspond to the code $3+2-1+4+5+$}
\end{figure}

\begin{Rem} On the figures we denote  the boundaries
of $\pi(H_i)$ by their numbers. The arrows point inward the
projections of the corresponding half-planes.
\end{Rem}

\subsubsection{Equivalent codes}

In the coding procedure described above we made several choices.
As a result we get several codes for the same situation. The
resulting classes are in fact orbits of a group acting on the set
of all possible codes.

This group is generated by two pairs of generators. The first pair
corresponds to the choices made on $ \mathbb{S}^1 $.

The first generator, denoted by $\b_1$, corresponds to the moving
the point $N$ to the previous  interval. It acts on the code by
cyclic permutation of the numbers and changing the sign of the
last element: the $i$-th number goes to the $(i+1)$-th place
except the first one which moves to the fifth place and changes
sign, e.g. $\b_1(1+2+3+4-5-)=5+1+2+3+4-$.

The second generator, denoted by $\b_2$, corresponds to the
change of orientations on the circle. It acts on codes by
symmetry: we should put the $i$-th number on the $5-i$-th place
preserving the sign (e.g. $\b_2(1+2+3+4-5-)=5-4-3+2+1+$).
\begin{figure}[h]
\centerline{\epsfysize=1in \epsffile{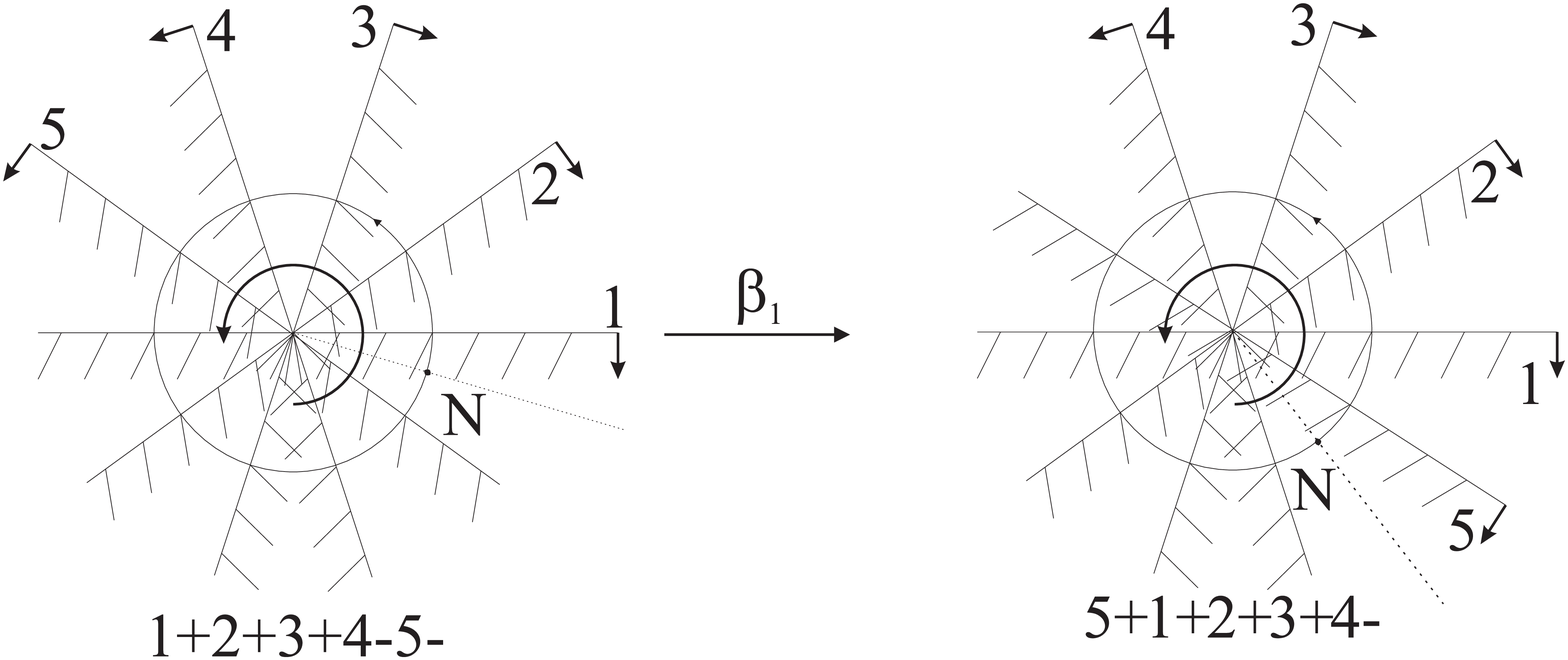}
\hfill\epsfysize=1in \epsffile{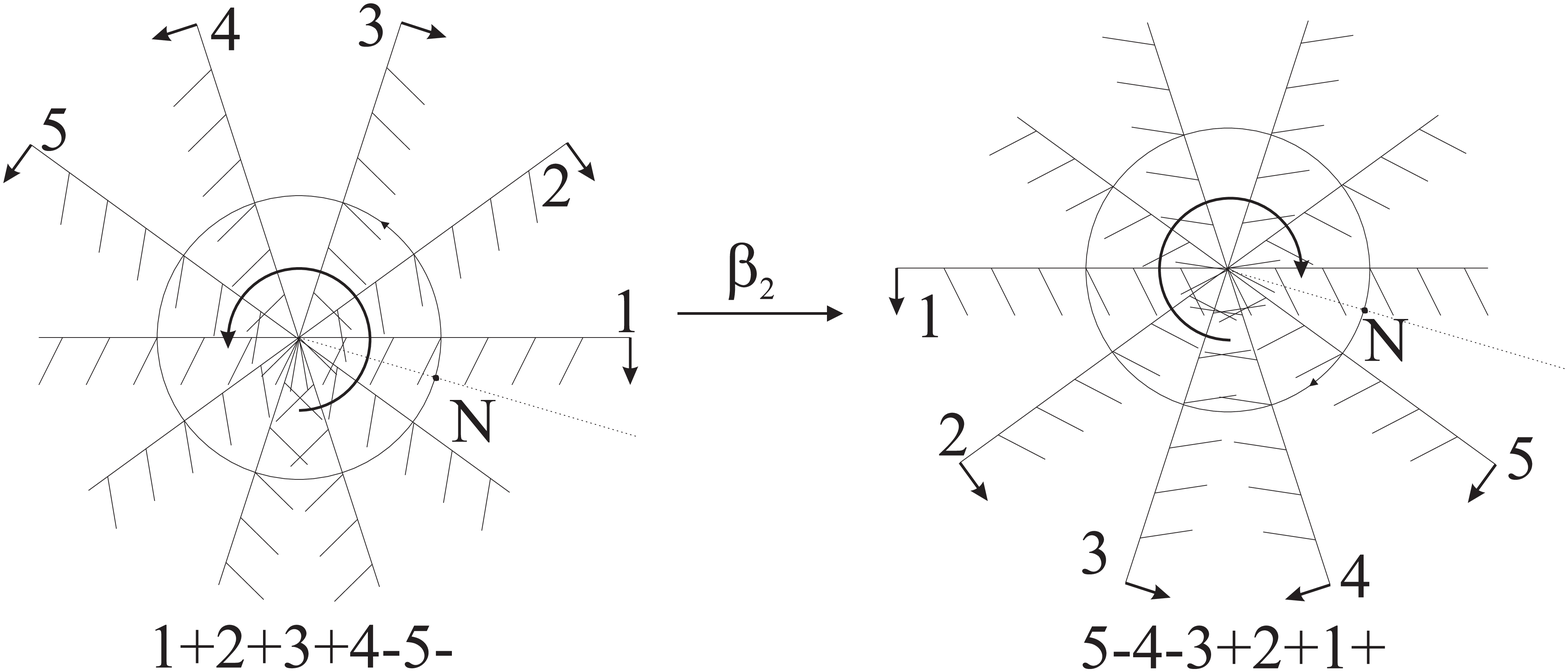}} \caption{Action of
$\b_1$ and $\b_2$}
\end{figure}

The second pair corresponds to the choices made on the Chebyshev
line $\ell$. In general, changing the position of the center of
the projection or the orientation results not only in change of
enumeration of half-planes but also  in the different choices of
the plane to which we project. So in order to describe the effect
of moving the point $M$ to the next interval or changing the
orientation of $\ell$ we have to identify somehow the planes of
projections.

The third generator of the group, denoted by $\a_1$, corresponds
to the moving the point $M$ to the point $M'$ in the previous
interval. If we identify planes $L_1$ and $L_2$ using the
projections  from $M'$ and make the same (upon this
identification)  choice of $N$ and of the orientation of $
\mathbb{S}^1$, then $a_1$ acts on codes by changing $1+$ to $2+$ ,
$2+$ to $3+$, $3+$ to $4+$, $4+$ to $5+$, $5+$ to $1-$, \dots,
$5-$ to $1+$ (e.g. $\a_1(1+2+3+4-5-=2+3+4+5-1+$). In other words,
the numeration shifts by $1$ and the image of the fifth (from the
$M$) half-plane flips.
\begin{figure}[h]
\centerline{\epsfysize=1.25in\epsffile{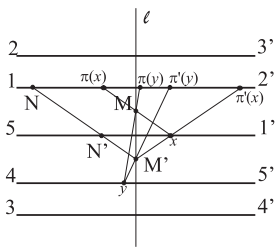}}
\caption{Projections from $M$ and $M'$ differ on $H_5$ and agree
on $H_1$, $H_2$, $H_3$ and $H_4$.}
\end{figure}

The fourth generator, denoted by $\a_2$, corresponds to the change
of orientation of $\ell$. After identifying $L_1$ and $L_5$ by
projection from $M$ action of $\a_2$ reduces to the renaming of
the planes. So $\a_2$ acts on codes by interchanging $5$ with $1$
and  $4$ with $2$ with  signs preserved (e.g.
$\a_2(1+2+3+4-5-=5+4+3+2-1-$).

\begin{figure}[h]
\centerline{\epsfysize=1in\epsffile{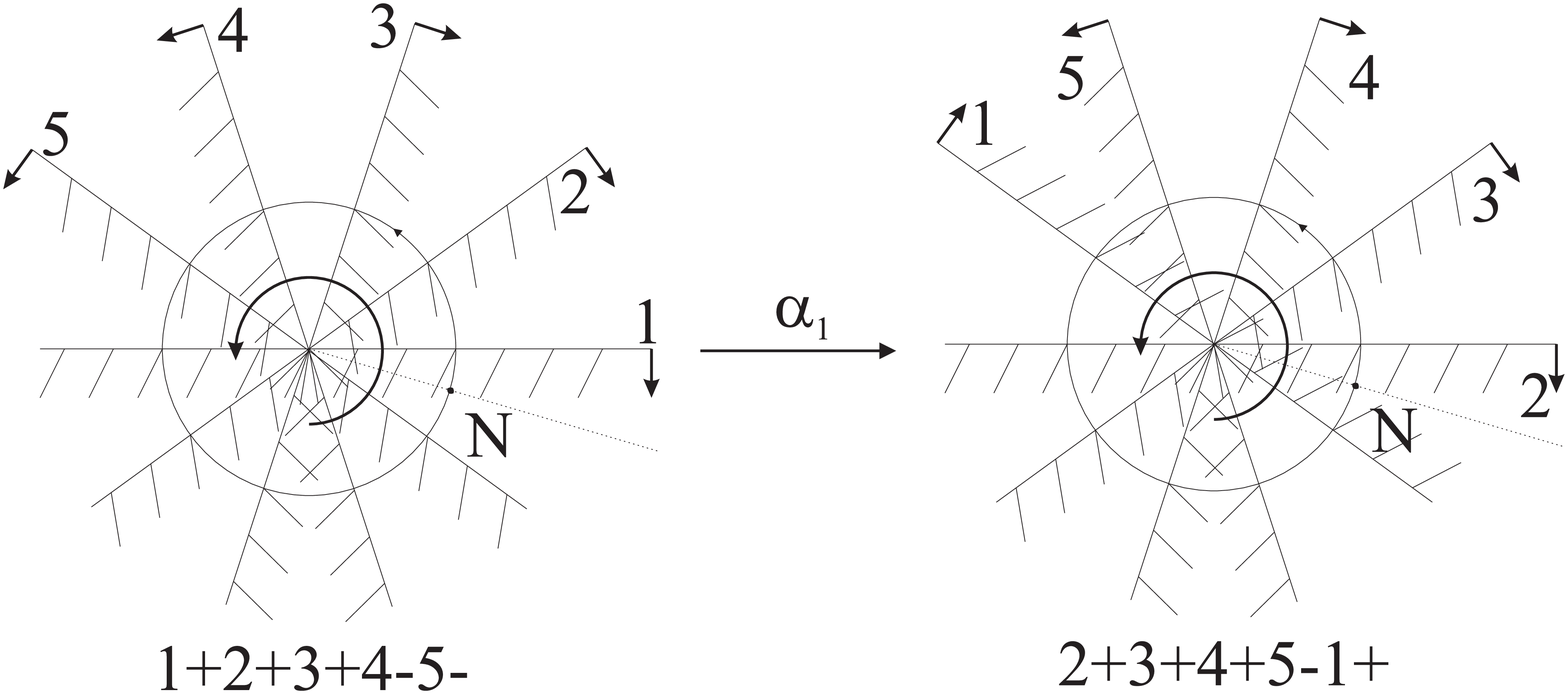}
\hfill\epsfysize=1in\epsffile{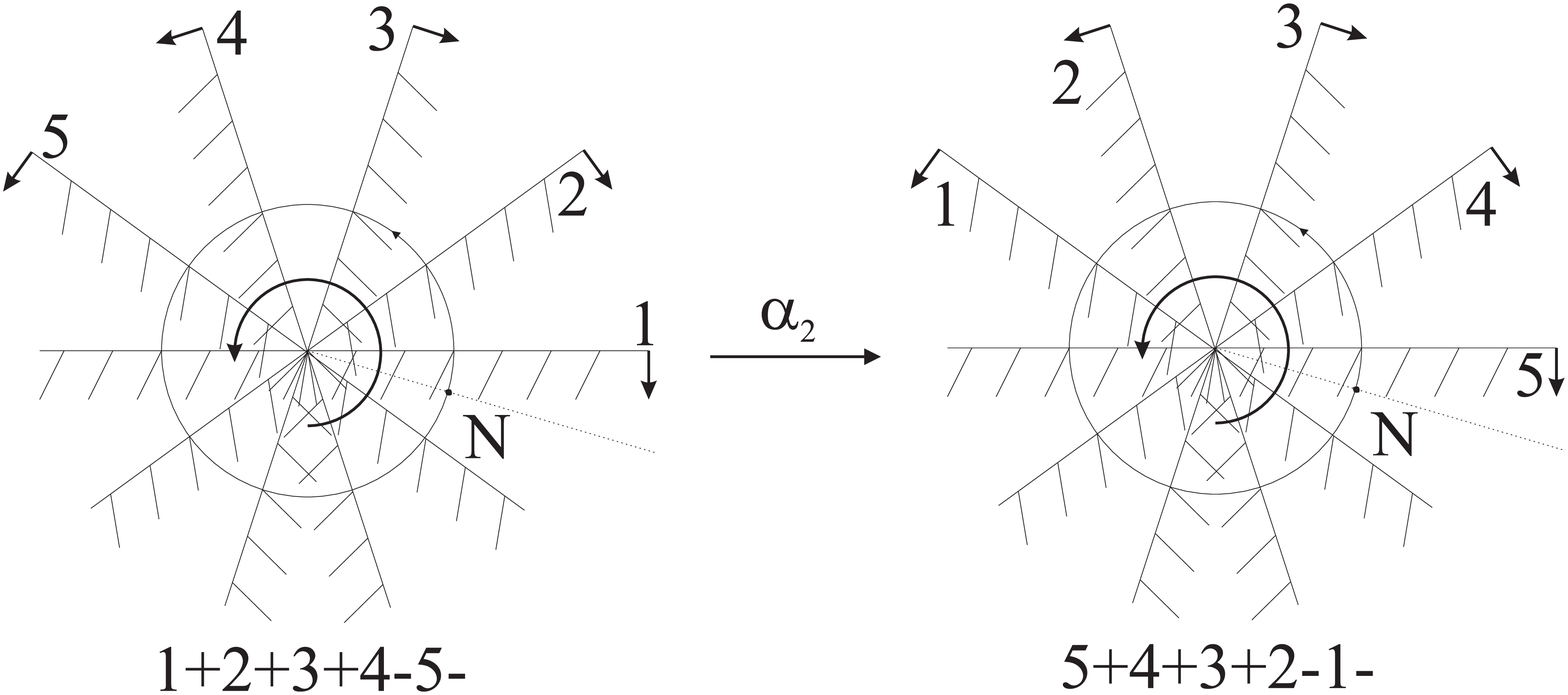}} \caption{Action of
$\a_1$ and $\a_2$}
\end{figure}

It is easy to see from this geometrical description that
$\alpha_i\beta_j\alpha_i^{-1}\beta_j{-1}=\alpha_1^{10}=\alpha_2^2=\beta_1^{10}=\beta_2^2=\operatorname{Id}$,
$\alpha_2\alpha_1\alpha_2=\alpha_1^{-1}$ and
$\beta_2\beta_1\beta_2=\beta_1^{-1}$, i.e. the group generated by
$\alpha_i$ and $\beta_i$ is $D_5\bigoplus D_5$.
\subsection{Cases of evident good deformation: trivial codes and Chebyshev property}
\subsubsection{Trivial codes}
There are cases (i.e. a combinatorial types of intersections of
projections of $H_i$) which are forbidden for Chebyshev lines.
These are in particular the cases when for some choice of $M$,
projections of all $H_i$ have nontrivial intersection (i.e. more
than one point). Indeed, in this case a good deformation which
will intersect the Chebyshev line can be easily found.

\begin{Thm}\label{triviality}
Configuration corresponding to a Chebyshev line cannot be coded by
a code containing $1+$, $2+$, $3+$ and $4+$.
\end{Thm}
\begin{proof}
 Suppose first that by choosing a point $M\in\ell$ and a point
$N\in\,\mathbb{S}^1$ we get a code consisting of  positive numbers
only, i.e a permutation of $1+$, $2+$, $3+$, $4+$ and $5+$. By
definition it means that the line connecting $M$ and $N$
intersects all $H_i$ at their interior, i.e. is a good
deformation.

If the code contains $5-$, then, after applying $\a_1^{-1}$, we
get an equivalent code with positive only entries, thus reducing
to the previous case.
\end{proof}

\subsubsection{Another easy case: the  Chebyshev
property}\label{ssec:Chebyshev property}

 The following lemma uses for the first (and the
last) time the Euclidean metric.  More exact, it uses the
definition of $H_i$ as the set of all points $x\in L_i$ such that
the scalar product
$(\overrightarrow{a_ix},\overrightarrow{a_is_i})$ is positive
(where $s_i\in S_i$ is the point of $S_i$ closest to $a_i$). We
will need this lemma only in the last chapter, when we consider
the six nontrivial codes.

\begin{Thm}\label{first in fourth}
No half-plane $H\subset L_1$ such that $a_1\in\partial H$ can
contain  $S_1$, $\pi(S_2)$, $\pi(S_3)$ and $\pi(S_4)$
simultaneously.
\end{Thm}
\begin{proof}
Denote by  $N$ the endpoint of inward normal $a_1N$ to $\partial
H$.

We are given  that $\pi(S_i)\subset H$ for $i=1,2,3,4$. Therefore
$\pi(s_i)\in H$, so, by \lemref{lem:euclid1}, $N\in \pi(H_i)$.

 If $N\not \in \pi(\partial H_i)$ then the
code corresponding to $N$ contains $1+$, $2+$, $3+$ and $4+$ and
we are done by the previous lemma.

If not, we can slightly move the point $N$ and get the same
result. Namely, suppose that $N\in \pi(\partial H_i)$ for some
$i$. Since (by genericity assumption) none of $\pi(\partial H_i)$
coincide, $N$ cannot lie on more than one $\pi(\partial H_i)$.
Therefore slightly moving $N$ inward this $\pi(H_i)$ we get a
point $N'$ corresponding to a code containing $1+$, $2+$, $3+$ and
$4+$, which is forbidden by the \thmref{triviality}.
\end{proof}

\begin{figure}[h]
\centerline{\epsfysize=0.25\vsize\epsffile{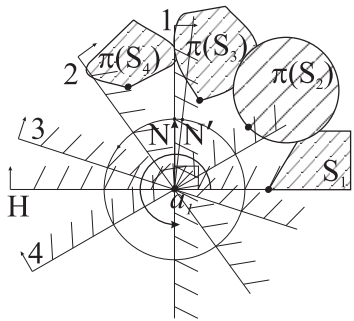}}
\caption{The case of $N\in\partial H_1:$ the point $N'$ lies in
the interiors  of $H_1, H_2, H_3$ and $H_4$.}
\end{figure}

\begin{Rem}
This lemma generalizes the following simple geometrical fact:
\begin{Lem}
There is no  half-space $H\subset \mathbb{R}^3$ with the Chebyshev
line on its boundary containing all five sections $S_i$.
\end{Lem}

Indeed, in this case in each plane $L_i$ we will get a figure like
in \lemref{lem:euclid1}, so a line obtained from a Chebyshev line
by a small parallel translation in the direction of the inward
normal to $\partial H$ will lie closer to all sections.
\end{Rem}

\section{Chess board}\label{boards}

In this section we single out all non-trivial codes, i.e. not
equivalent to the named  in \thmref{triviality}. Though the number
of codes is huge (namely $3840=2^55!$), there are only six
equivalency classes not containing  trivial  codes. They are
listed in the \thmref{nontrivial codes} below.

\subsection{From a code to a corresponding chessboard}
It is easier to visualize codes as a position of five rooks on a
$5\times5$ chess board. This is done as follows: in the first
column we put the rook in the row which number is equal to the
first number in the code. The color of the rook is white if this
first number has sign $+$ and black otherwise. We continue like
this for the second, third, fourth and fifth column (so if we
forget the colors, the rooks position is exactly the graph of the
permutation given by the code).

\begin{figure}[h]
\centerline{\epsfysize=0.2\vsize\epsffile{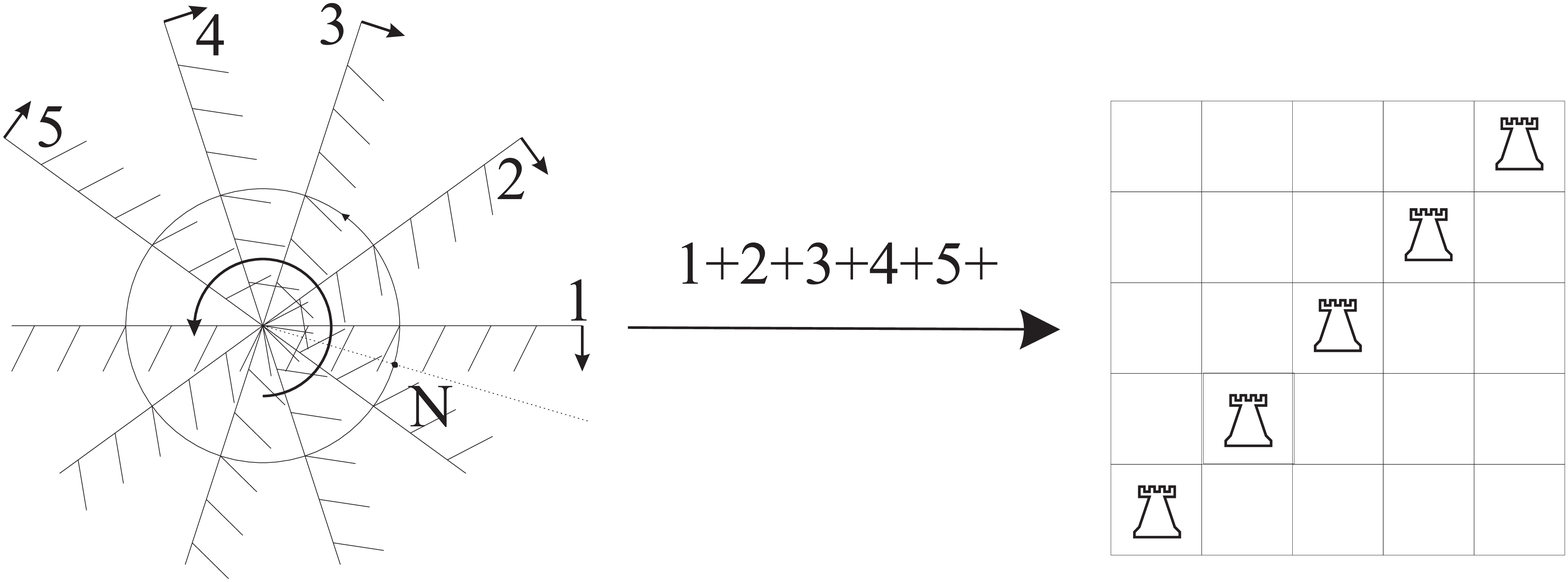}}
\caption{The code and the board corresponding to the projection
above.}
\end{figure}

 It is easy to see that each column or row contains
exactly one rook, i.e. the rooks do not threaten  each other.

\subsection{How the symmetry group acts on rooks positions}
We described above  an action of some symmetry group on codes. In
the chess board realization the action of this  group is
remarkably simple:
\begin{itemize}
 \item $\b_1$ acts by moving  the fifth column to the first
place and
 changes the color of the rook standing in this column;
 \item $\a_1$ acts in a similar way but with rows: $\a_1$ moves
the fifth row
 to the first place  and
 changes  the color of the rook standing in this row;
 \item $\a_2$ acts by symmetry with respect to the vertical
line;
 \item $\b_2$ acts by symmetry with respect to the horizontal
line.
\end{itemize}
\begin{figure}[h]
\centerline{\hfill\epsfysize=1.2in
\epsffile{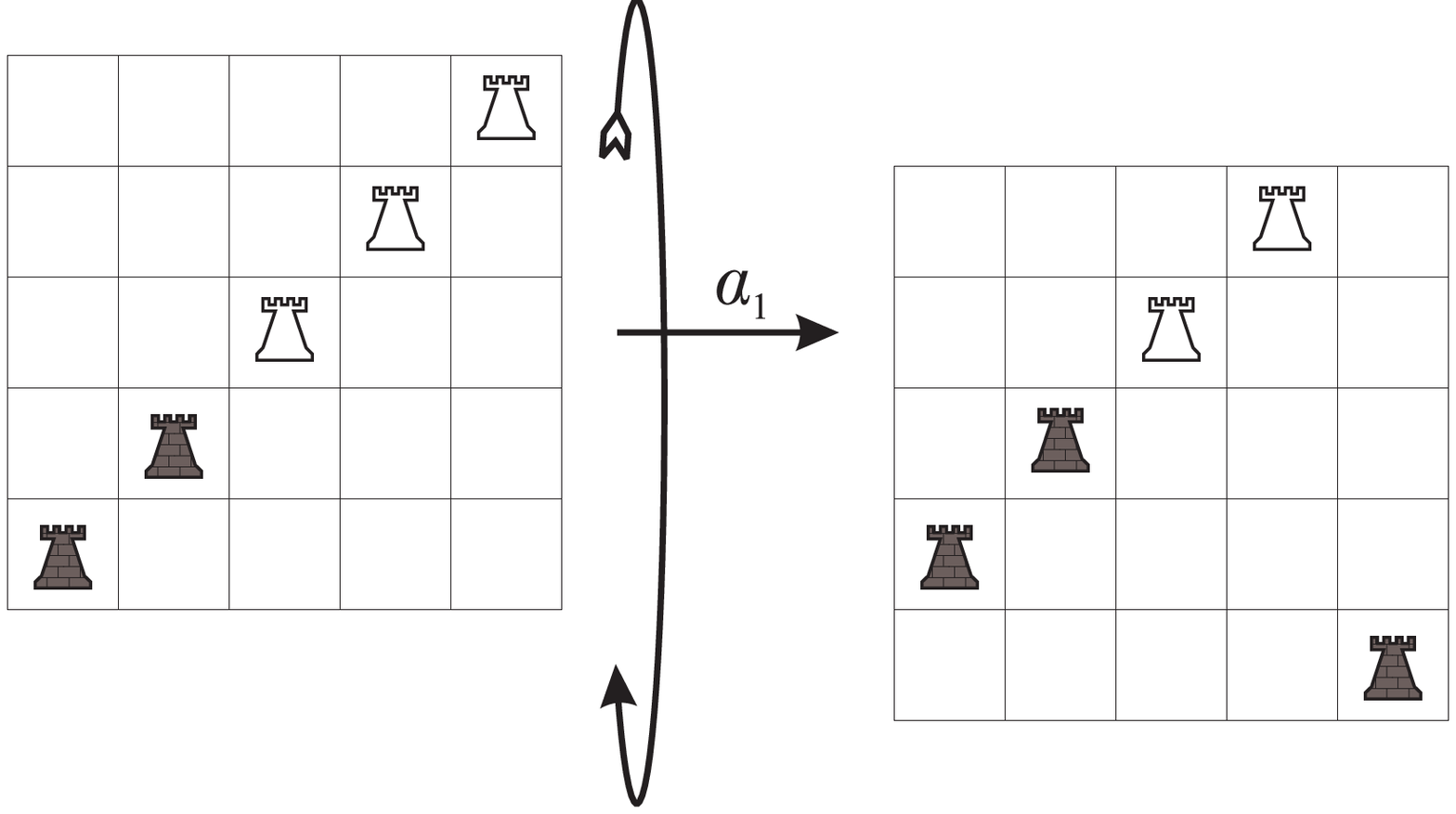}\hfill\epsfysize=1.2in\epsffile{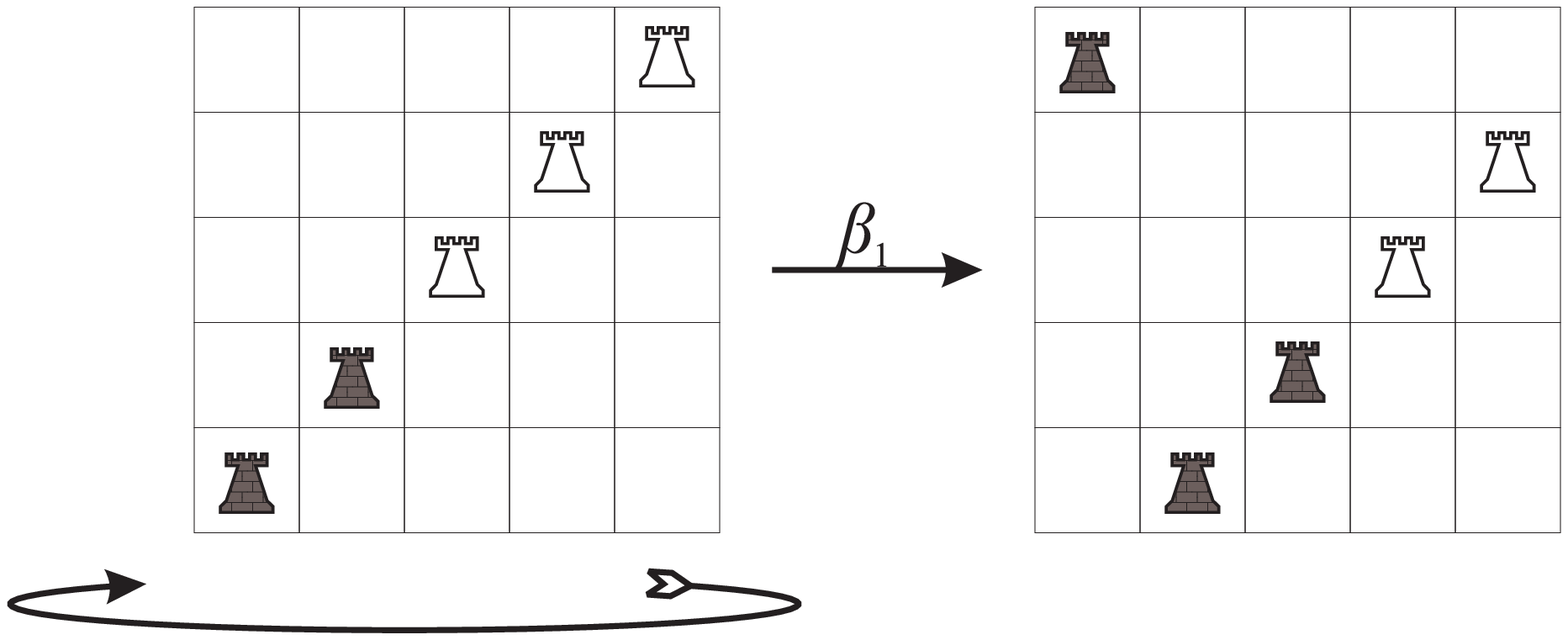}\hfill}
\vskip 0.5cm
\centerline{\hfill\epsfysize=1.2in\epsffile{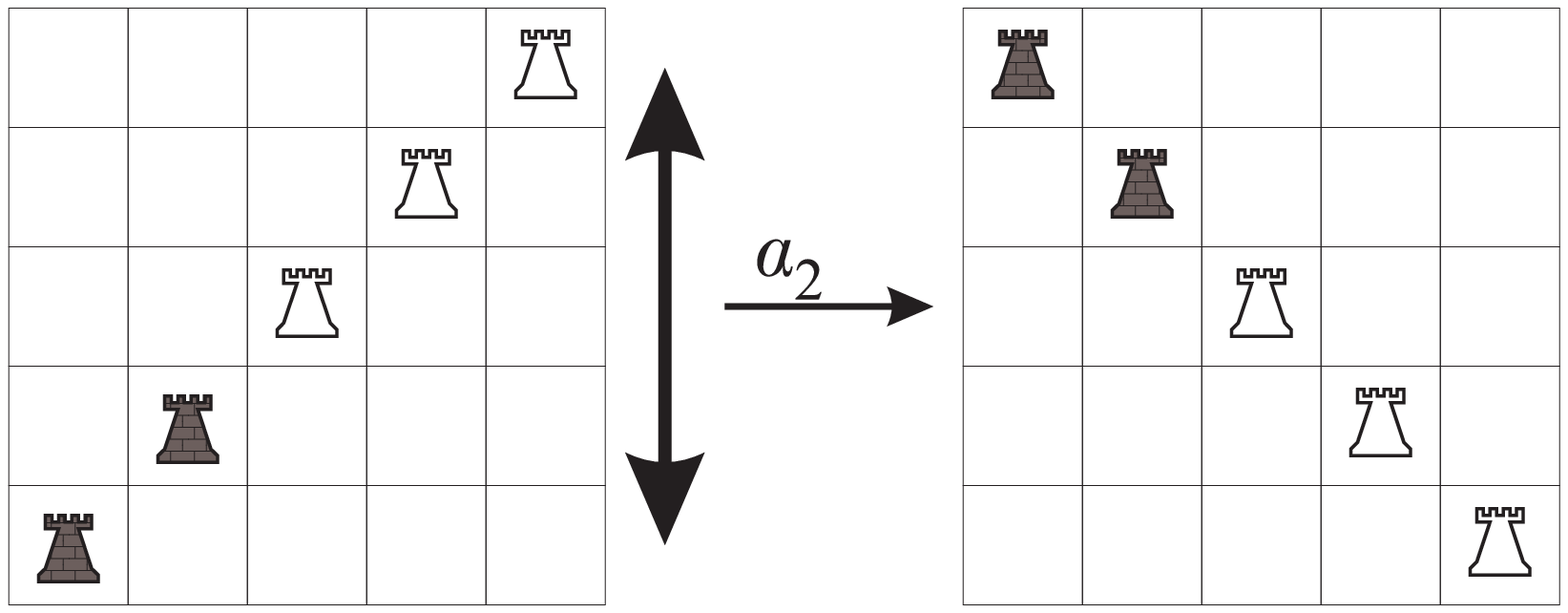}\hfill\epsfysize=1.2in\epsffile{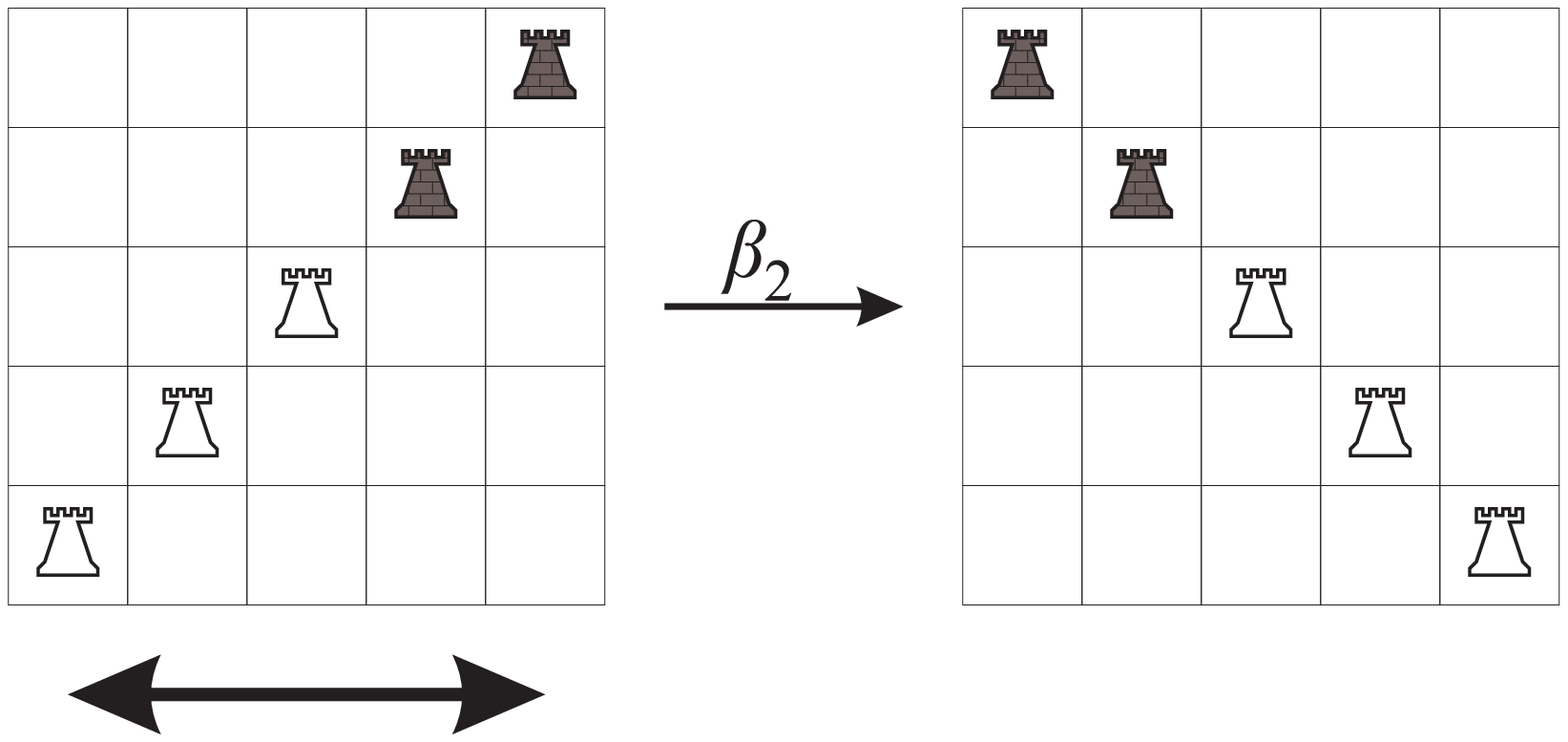}\hfill}
\caption{Action of generators of the group}
\end{figure}

\subsection{Six equivalence classes consisting of nontrivial arrangements only}

The trivial codes correspond to the arrangements of white rooks
only, which will be called trivial arrangements. Our goal is to
exclude rooks arrangements equivalent to trivial ones. This is
done in this subsection by a

\begin{Lem} Any arrangement non-equivalent to a trivial case
is equivalent to a arrangement with only one black rook.
Moreover,
this rook can be supposed to stand not on the border of the
board.
\end{Lem}
\begin{proof}
Pick any arrangement which is not equivalent to a trivial one. The
$\b_1^5$ simply changes all colors to the opposite ones, so we can
assume that the number of black rooks is equal to one or two. The
first case is what we need, so suppose that there are two black
rooks. If one of them stands on the first or the last row, then
using $\a_1^{\pm 1}$ we can change its color without changing the
color of others, so leaving only one black rook. Similar statement
holds for columns and $\b_1$.

So we can suppose that both black rooks are in the inner
$3\times3$ square. Then we get at least two white rooks on the
border. Take the fifth row. It contains one rook. Therefore a
first or a fifth column should contain another white rook and
moving this column and the fifth row (i.e. acting by $\b_1\a_1$ or
by $\b_1^{-1}\a_1$) we arrive to a situation with four black
rooks, which is equivalent (by $\b_1^5$) to a situation with one
rook only.

This black rook cannot stand on the border since otherwise by one
move ($\a_1^{\pm 1}$ or $\b_1^{\pm 1}$) we arrive to a trivial
situation.
\end{proof}

Using the symmetries $\a_2$ and $\b_2$, we can  assume that the
black rook occupies one of the four squares
$(2,2),(2,3),(3,2),(3,3)$.

\subsubsection{The case $(2,2)$}
Consider first the case of  the black rook on the square
$(2,2)$.
\begin{Lem}\label{lem:2,2 and corners}
If one of the squares $(1,1),(1,5),(5,1)$ is occupied, then the
position is trivial.
\end{Lem}
\begin{proof}
Indeed, in these cases $\b_1^{-1}\a_1^{-2}$ or
$\b_1^{-2}\a_1^{-1}$ or  $\b_1\a_1^{-2}$ correspondingly
transforms the position to a trivial one.
\end{proof}

Therefore in a position not equivalent to a trivial one the white
rook of the first column can occupy one of the squares $(1,3)$ or
$(1,4)$ only and the square $(5,1)$ is empty.

\subsubsection{White on $(1,3)$ and Black on $(2,2)$}

This leaves four configurations:
\begin{itemize}
\item[C1] $3+2-1+4+5+$
\item[C2] $3+2-1+5+4+$
\item[C3] $3+2-4+1+5+$
\item[C4] $3+2-5+1+4+$
\end{itemize}
\vskip 0.5cm
\begin{figure}[h]
\centerline{\hfill\epsfysize=0.11\vsize\epsffile{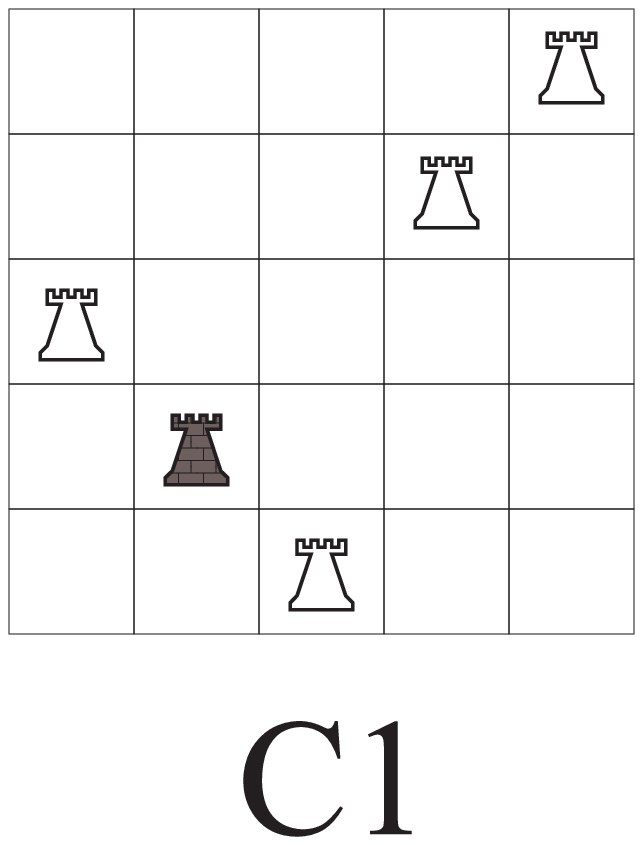}\hfill\epsfysize=0.11\vsize\epsffile{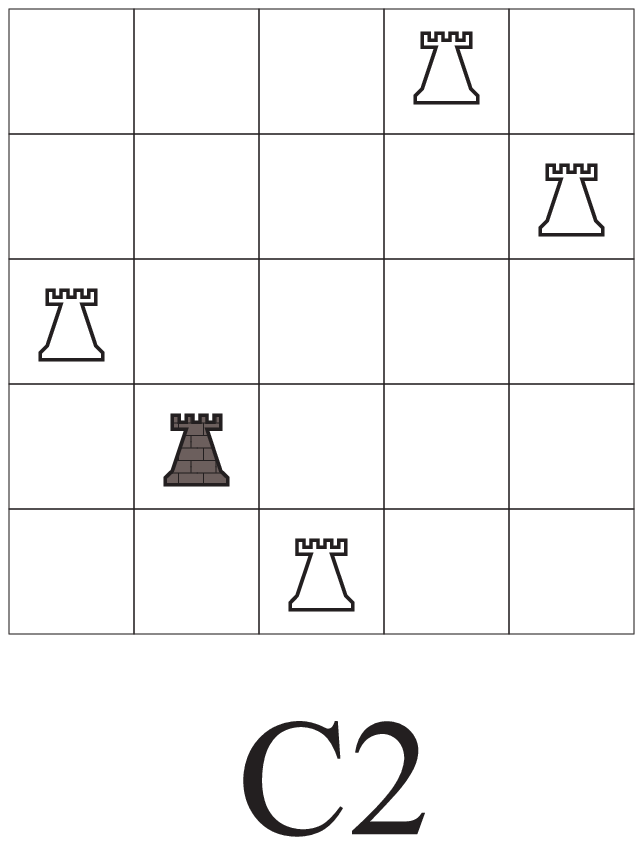}
\hfill\epsfysize=0.11\vsize\epsffile{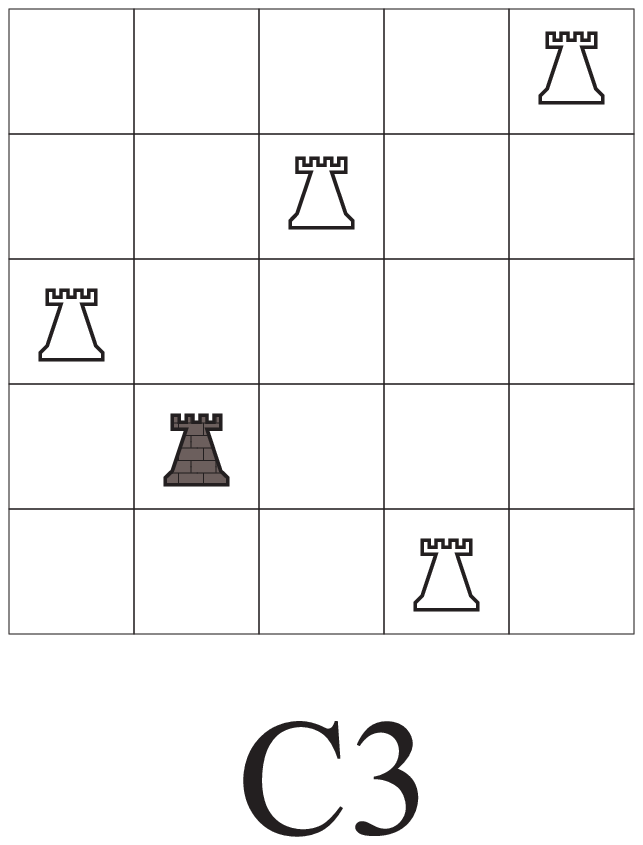}\hfill\epsfysize=0.11\vsize\epsffile{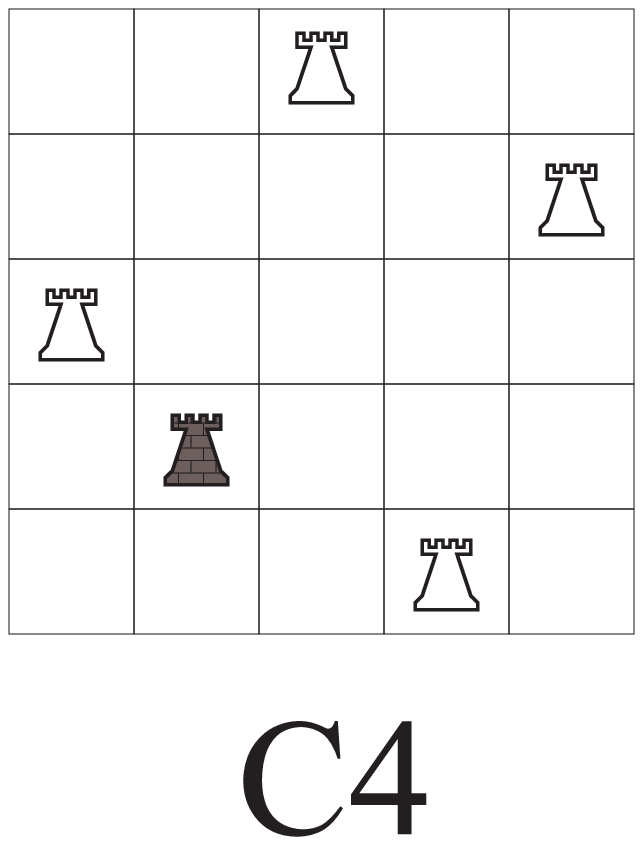}\hfill}
\end{figure}

\subsubsection{White on $(1,4)$ and Black on $(2,2)$}

These are another four possibilities (remind that $(5,1)$ is
empty):
\begin{itemize}
\item[D2] $4+2-1+5+3+$
\item[D3] $4+2-3+1+5+$
\item[D4] $4+2-5+1+3+$
\item[D5] $4+2-1+3+5+$
\end{itemize}

But D4 becomes trivial after $\a_1^{-3}\b_1^2$, and D2 becomes  D3
after $\b_1^2\a_2\b_2$. Moreover, after $\a_1^3\b_1^{-3}$ C4
becomes D2. So the only new configuration is D5.
\begin{figure}[h]
\centerline{\epsfysize=1in \epsffile{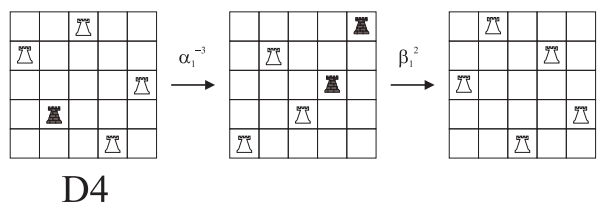}} \caption{ D4 becomes
trivial after $\a_1^{-3}\b_1^2$.}
\end{figure}

\begin{figure}[h]
\centerline{\epsfysize=0.1\vsize\epsffile{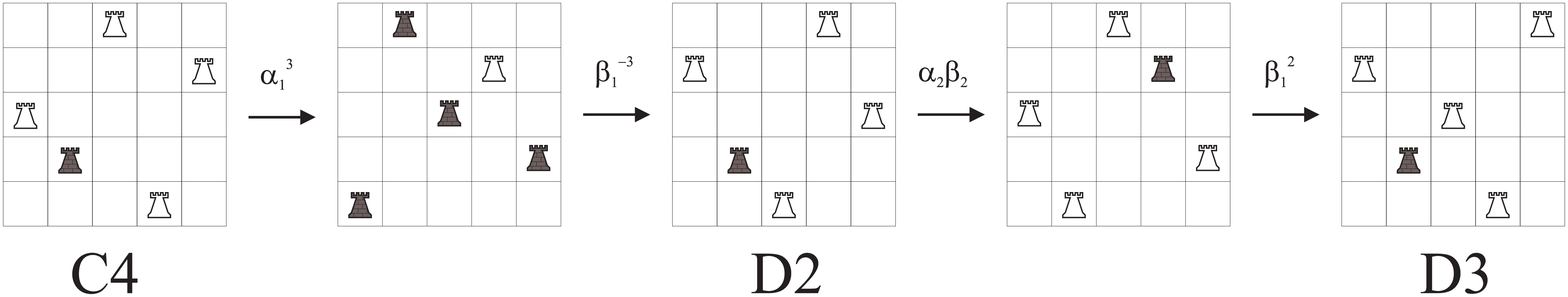}}
\caption{ C2 is equivalent to D2 and to D3.}
\end{figure}

\subsubsection{Black on $(2,3)$}

Similarly to \lemref{lem:2,2 and corners}, the white rook in the
first column cannot stand on the first or the last row. In other
words, in a position with black rook on $(2,3)$ and not equivalent
to a trivial one the squares $(1,5)$ and $(1,1)$ are empty.
Indeed, $\a_1^{\pm1}\b_1^{-2}$ correspondingly trivialize these
arrangements.

So the only places the  white rook can stand on are $(1,2)$ or
$(1,4)$. These positions are in fact equivalent by $\a_2$, so we
can consider the positions with a white rook on $(1,4)$ and the
black rook on $(2,3)$.

But these positions are equivalent by $\a_2\b_2\b_1^{-2}$ to the
positions with the black rook on $(2,2)$, so are in fact
considered above.

\begin{figure}[h]
\centerline{\epsfysize=0.13\vsize\epsffile{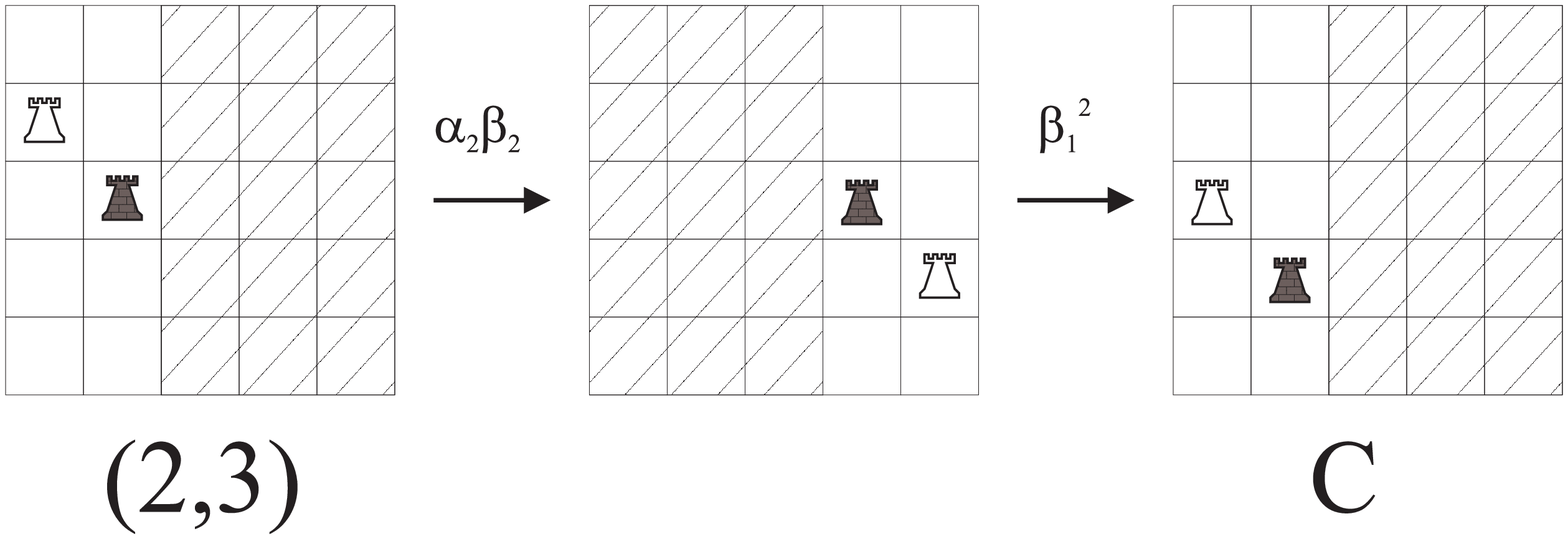}} \caption{
The series of cases of the black rook on $(2,3)$ is equivalent to
the series C.}
\end{figure}

\subsubsection{Black on $(3,2)$}
These arrangements are also equivalent to arrangements with the
(only) black rook on $(2,2)$. The proof repeats word-by-word the
proof above with change of $\b$ to $\a$ and of $\a$ to $\b$
everywhere. This is because the actions of the group is symmetric
with respect to diagonal (though this symmetry isn't  itself in
the group).
\subsubsection{Black rook on $(3,3)$}
The complement of the square to the third row and the third column
consists of four two-by-two squares.
\begin{Lem}
If the arrangement is not equivalent to a trivial one, then each
square contains exactly one rook.
\end{Lem}
\begin{proof}
Indeed, if not, then one of them contains two rooks and the
opposite should necessarily contain the other two (since in each
row and in each column stands exactly one rook). Applying $\b_2$
if necessary, one can suppose that these are the lower left and
the upper right squares. Then $\a_1^2\b_1^3$ transforms
arrangement to a trivial one.
\end{proof}

\begin{figure}[h]
\centerline{\epsfysize=0.13\vsize\epsffile{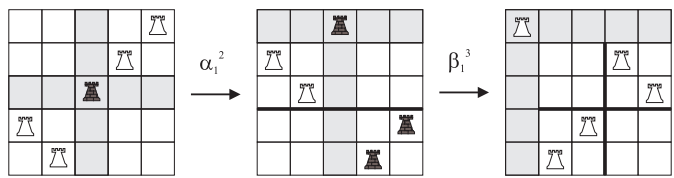}}
\caption{ Triviality of the case of the black rook on $(3,3)$ and
one of the squares containing two rooks.}
\end{figure}
\begin{Lem}
If one of rooks stands in the corner  (i.e. on $(1,1),
(1,5),(5,1)$ or $(5,5)$),  then the situation is equivalent to a
situation with the only black rook standing on $(2,2)$ (i.e. is in
fact considered above).
\end{Lem}
\begin{proof}
Using $\a_2$ and $\b_2$, if necessary, we can suppose that the
white rook stands on $(1,1)$. Then we get a situation with the
only black rook on $(2,2)$ after $\a_1^{-1}\b_1^{-1}$.
\end{proof}
\begin{Cor}
All configurations with one of  white rooks in the inner $3\times
3$ square are trivial or have a rook in a corner.
\end{Cor}
\begin{proof}
Suppose that $(2,2)$ is occupied and the position is neither
trivial nor with a rook in a corner. Then the $2\times 2$ square
contain one rook each. Then the squares $(1,4)$ and $(4,1)$ are
occupied, since the corners are empty and the second row and
second column already contain a rook. Therefore the only remaining
square for the fourth rook is in the corner $(5,5)$, which is
forbidden.\end{proof}
\begin{figure}[h]
\centerline{\epsfysize=0.13\vsize\epsffile{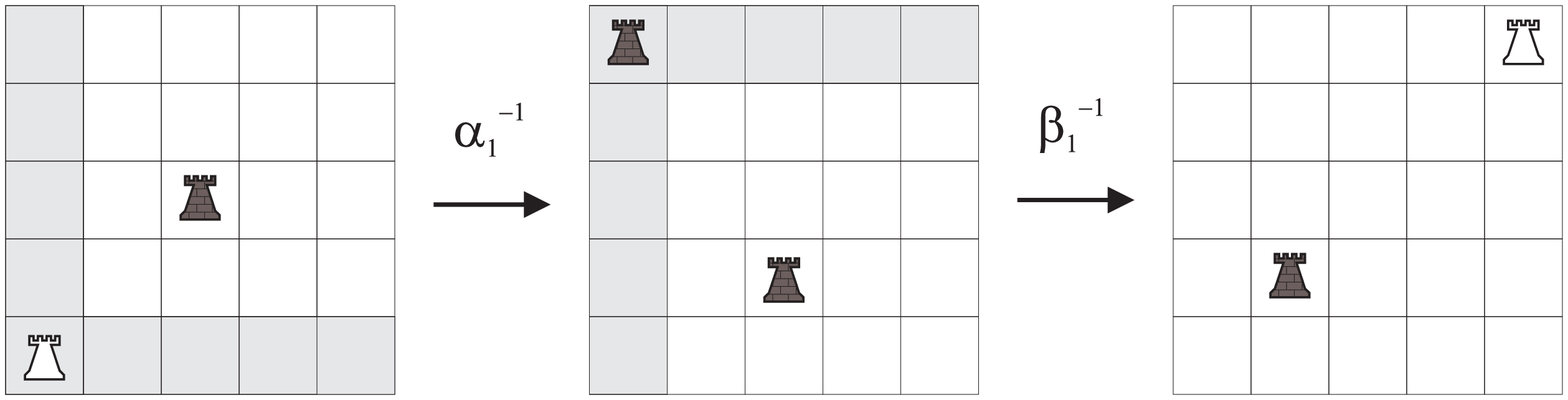}}\caption{
The case of the black rook on $(3,3)$ and one of the rooks in the
corner is equivalent to the series C.}
\end{figure}

The only remaining positions are $4+1+3-5+2+$ and $2+5+3-1+4+$,
which are equivalent by $\a_2$ or $\b_2$.

\subsubsection{The final list}
It consists of six variants.

\begin{Thm}\label{nontrivial codes}
A configuration corresponding to a Chebyshev line should  be
equivalent to a configuration described by one of the following
codes

\begin{tabular}{ccc}
C1 $3+2-1+4+5+$ & C2 $3+2-1+5+4+$& C3 $3+2-4+1+5+$\\
C4 $3+2-5+1+4+$ & D5 $4+2-1+3+5+$ &E6 $4+1+3-5+2+$
\end{tabular}
\end{Thm}

\begin{figure}[h]
\centerline{\epsfysize=0.22\vsize {\epsffile{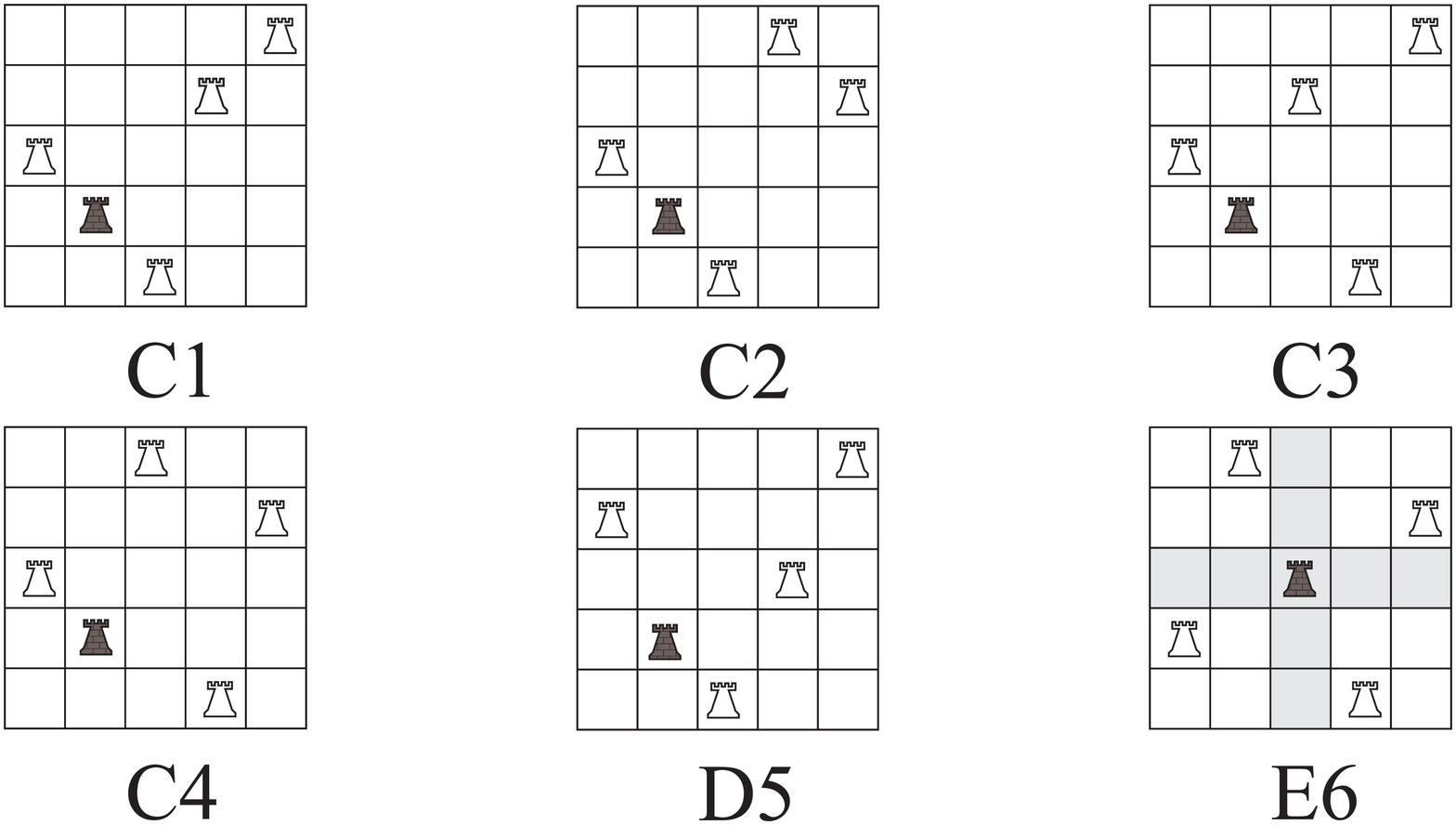}}}
\caption{The six nontrivial variants}
\end{figure}

\section{Non-triviality of a code and convex-concavity imply existence of good deformation}\label{final}
 In this chapter we consider the six
nontrivial cases of \thmref{nontrivial codes}. Each case has
several continuous parameters (e. g. angles between $\partial
H_i$, distances between $L_i$), and only  for some  choice of
parameters the configuration of half-planes arises from a
Chebyshev line. In other words, for only part of the parameter
space parameterizing this combinatorial type the corresponding
configuration of half-planes do not admit a good deformation.
Indeed, the \thmref{nontrivial codes} excludes only codes
admitting a good deformation {\em intersecting the Chebyshev line
$\ell$\/}, and do not deal with good deformations not intersecting
$\ell$.

In what follows we show that the configurations of half-planes
arising from   sections $S_i$ of a convex-concave body all admit a
good deformation. Therefore they cannot correspond to a Chebyshev
line, so the assumption that the Chebyshev line doesn't intersect
the sections leads to a contradiction.

More exact, we extract from the convex-concavity condition some
inequality between double ratio of angles between $\partial H_i$
and double ratios of distances between $L_i$ in some particular
combinatorial assumptions. This inequality implies existence of a
line intersecting four from half-planes $H_i$ in some particular
sectors. For five from the six cases of \thmref{nontrivial codes}
these assumptions are satisfied, and moreover the resulting line
automatically intersects the fifth half-plane. The sixth case  E6
simply cannot occur for convex-concave sections.

The main tool in the proofs is the \thmref{Browder}, only applied
now to some parts of the sections $S_i$. The only Euclidean
property we will need is the \thmref{first in fourth}, which
statement is projective. So we can move  the center of projection
to infinity, and  the projection becomes  a parallel projection
$\pi:\mathbb{R}^3 \to L_1$ along the $z$-axis, with $S_i$ are
ordered by their $z$-coordinate.

We will also use a linear structure defined on $L_1$ defined by
the coordinates $x$ and $y$ (i.e. we take the point $a_1$ as the
origin).

\subsection{Sectorial Browder Theorem}

We will denote by $\pi(H_i)^c$ for the  closure of
$L_1\setminus\pi(H_i)$. We define half-spaces
$B_i=\pi^{-1}(\pi(H_i))$ and denote by $B_i^c$ the closure of
their complements.
\begin{Thm}\label{Restricted Brauer}
 Suppose that
\begin{enumerate}
\item $H_1\cap \pi(H_4)^c\subset \pi(H_2)^c$ and
\item $\pi(H_3)\cap \pi(H_2)\subset \pi(H_4)^c$.
\end{enumerate}
Suppose moreover that $S_1\cap \pi(H_4)^c\not=\emptyset$. Then
$\pi(S_2)\cap \pi(H_3)$, $\pi(S_3)\cap \pi(H_2)$ and
$\pi(S_4)\cap
H_1^c$ are also non-empty and there exists a straight line $L$
intersecting $S_1\cap B_4^c$, $S_2\cap B_3$, $S_3\cap B_2$ and
$S_4\cap B_1^c$.
\end{Thm}
\begin{figure}[h]
\centerline{\epsfysize=2in
\epsffile{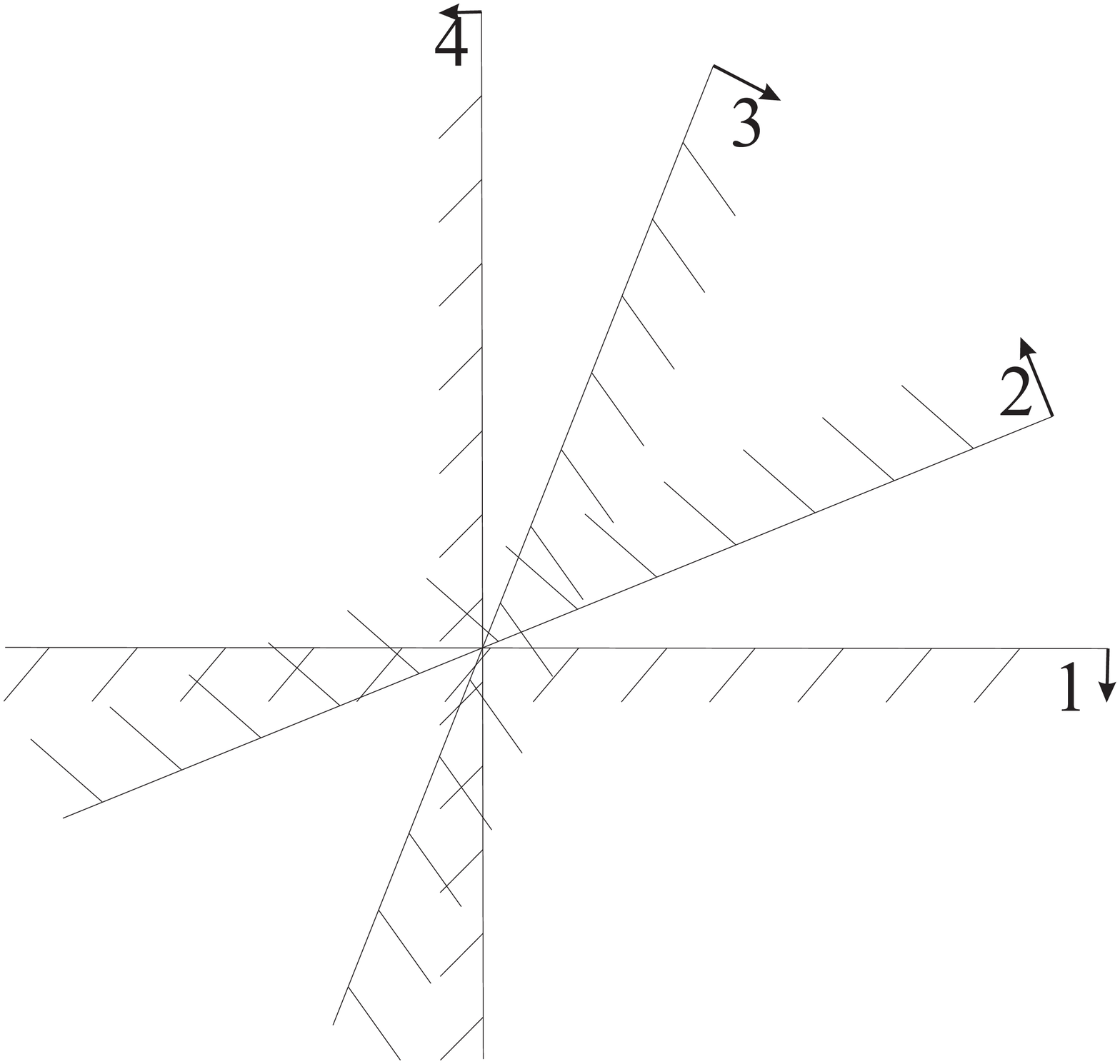}\hfill\epsfysize=2in
\epsffile{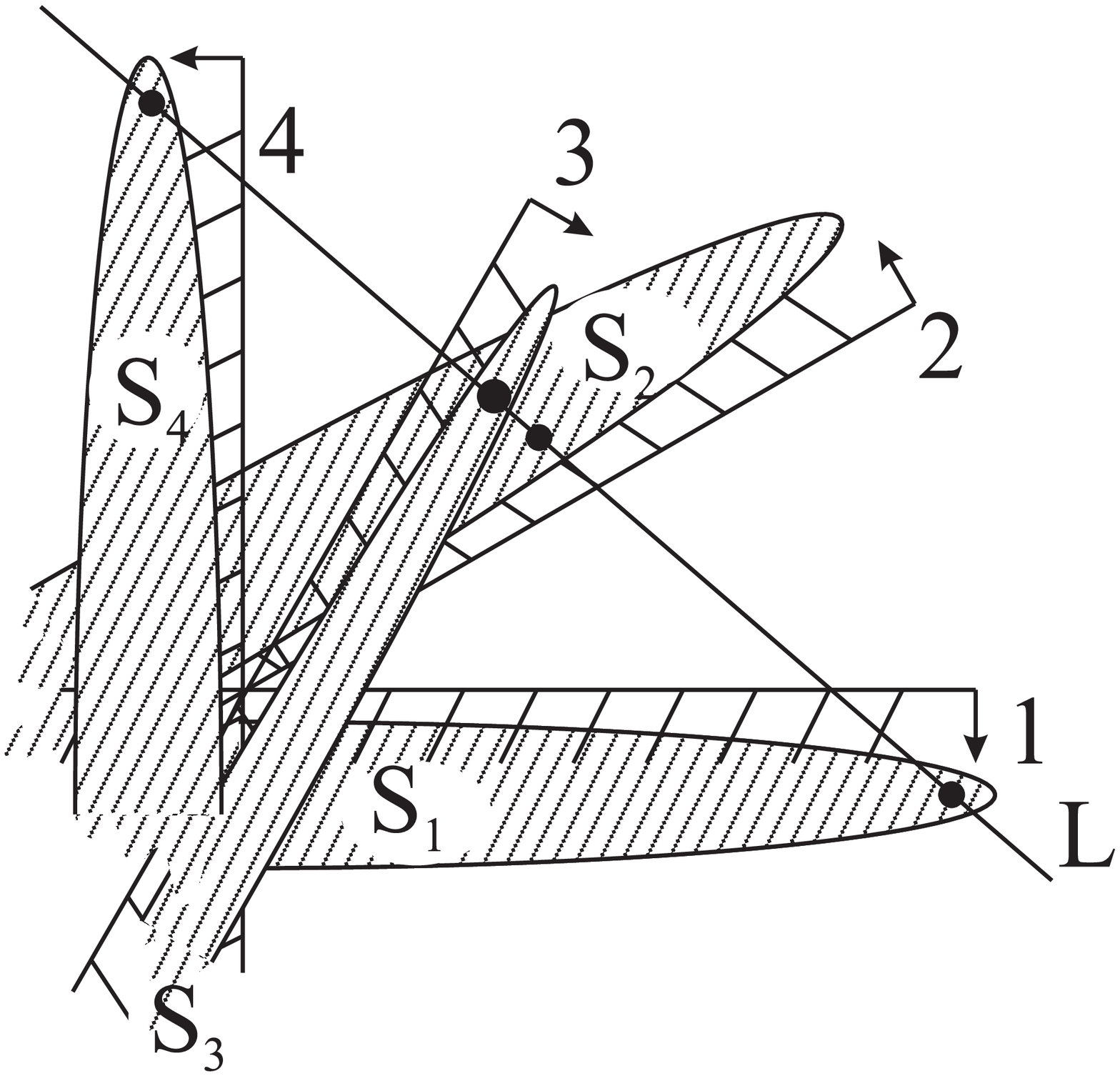}} \caption{The configuration of
half-planes and projection of the line $L$ from the
\thmref{Restricted Brauer}. }
\end{figure}

In our notations the conditions (1) and (2) mean existence of the
subsequence $1+2-3+4-$ in a sequence coding the configuration. In
applications below the condition $S_1\cap
\pi(H_4)^c\not=\emptyset$ will follow from the \lemref{trivial-not
trivial} below.
\begin{proof}
First we prove two combinatorial lemmas:
\begin{Lem}
$H_1\cap \pi(H_4)^c\subset \pi(H_3)$
\end{Lem}
\begin{proof}
Suppose that $H_1\cap \pi(H_4)^c\not\subset \pi(H_3)$. Since
boundaries of the half-planes are pairwise different, there is a
point $x$ lying in the interior of $(H_1\cap \pi(H_4)^c)\setminus
\pi(H_3)$. Then $-x\in H_1^c\cap \pi(H_4)\cap \pi(H_3)\subset
\pi(H_2)\cap \pi(H_3)$ by assumption and also $-x\in \pi(H_4)$ -
contradiction.
\end{proof}
\begin{Lem}
$\pi(H_2)\cap \pi(H_3)\subset \pi(H_1)^c$
\end{Lem}
\begin{proof}
As before, take $x$ in the interior of $\pi(H_2)\cap \pi(H_3)\cap
H_1$. Then $x\in \pi(H_4)^c\cap \pi(H_1)$ by the  assumption (2)
and therefore $x\in \pi(H_2)^c$ by the  assumption (1) -
contradiction.
\end{proof}

Our claim will be proved by applying the \thmref{Fan} to $S_1\cap
\pi(H_4)^c$ as $B$, $S_2\cap B_3$ as $A$, $S_3\cap B_2$ as $C$ and
$S_4\cap B_1^c$ as $D$. Let's  check conditions of \thmref{Fan}.
In other words, we have to check that
\begin{enumerate}
\item a line passing through $S_1\cap B_4^c$ and intersecting
$S_2$ and $S_3$ (existing by convex-concavity) intersects $S_2\cap
B_3$ and $S_3\cap B_2$ and
\item a line passing through $S_3\cap B_2$ and intersecting
$S_1$ and $S_4$ (existing by convex-concavity) intersects $S_4\cap
B_1^c$ and $S_1\cap B_4^c$.
\end{enumerate}
(Clearly $S_1\cap B_4^c$, $S_2\cap B_3$, $S_3\cap B_2$ and
$S_4\cap B_1^c$ are compact and convex).

Let a line intersects $S_1\cap B_4^c$ and $S_2$ and $S_3$ at
points $c_1, c_2$ and $c_3$ accordingly. Necessarily $c_2$ lies
between $c_1$ and $c_3$. We know that $c_1\in S_1\cap
\pi(H_4)^c\subset H_1\cap \pi(H_4)^c\subset \pi(H_2)^c\cap
\pi(H_3)$. Since $c_1,c_3\in B_3$, so $c_2\in B_3$ (so $S_2\cap
B_3$ is non-empty). Similarly, $c_1\in B_2^c$ and $c_2\in B_2$, so
$c_3\in B_2$ (and $S_3\cap B_2$ is non-empty). So the first claim
is proved.

Similarly, let a line intersects $S_3\cap B_2$ and $S_4$ and $S_1$
at points $c_3, c_4$ and $c_1$ accordingly. As before, $s_3\subset
B_4^c\cap B_1^c$.  Since $c_4\in B_4$ and $c_3\in B_4^c$, so
$c_1\in S_1\cap \pi(H_4)^c$. Since $c_3\in B_1^c$, so $c_4\in
S_4\cap B_1^c$ (so in particular $S_4\cap B_1^c$ is not empty).
The second claim follows.

\end{proof}

\subsection{Double ratios}

After projecting a configuration satisfying conditions of the
\thmref{Restricted Brauer} to the plane $L_1$ we obtain a figure
below.
\begin{figure}[h]
\centerline{\epsfysize=1.5in\epsffile{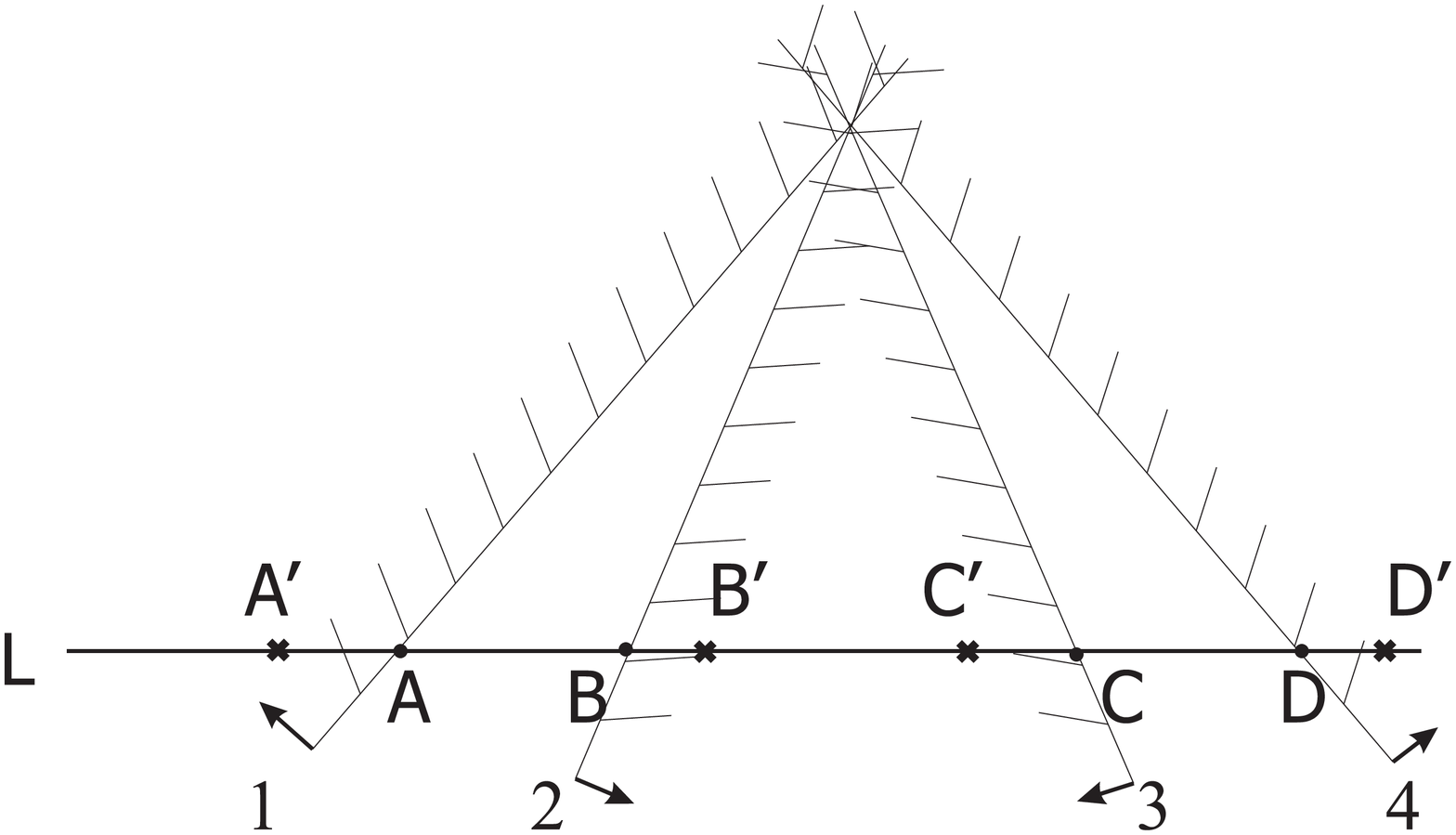}}
\end{figure}

 Here $L$ is the projection of the line existing
by \thmref{Restricted Brauer}. By $A', B', C', D'$ we denote
intersections of $L$ with $S_i$ and by $A,B,C,D$ intersections of
$L$ with $\partial (\pi(H_i))$.

The existence of the line $L$ implies an inequality between the
double ratio of distances between $L_i$ and the double ratio of
directions of boundaries of $H_i$. Namely, denote the double ratio
$\frac{AB}{BD}\colon\frac{AC}{CD}$ of points $A,B,C,D$ by
$(A,B,C,D)$. Then $(A'B'C'D')$ is exactly the double ratio of
distances between $L_i$:
$$
(A,B,C,D)=\frac{h_1}{h_1+h_2}\colon\frac{h_1+h_2}{h_3},
$$
where $h_i$ is the distance between $L_i$ and $L_{i+1}$.
$(A,B,C,D)$ is the double ratio of directions of $\partial H_i$
and the following inequality holds:
\begin{Cor}\label{double-ratios relation}
In assumptions of \thmref{Restricted Brauer} the double ratio of
distances between $L_i$ is strictly smaller than the double ratio
of directions of $\partial H_i$:
$$ (A',B',C',D')>(A,B,C,D)
$$
\end{Cor}
\begin{proof}
Indeed, the configuration of the points $A',B',C',D'$ is obtained
from the points $ABCD$ by the movements which only increase the
above double ratio:
\begin{enumerate}
\item $(A,B,C,D)<(A',B,C,D)$ since
$\frac{A'B}{A'C}>\frac{AB}{AC}$,
\item $(A',B,C,D)<(A',B',C,D)$ since
$\frac{A'B'}{B'D}>\frac{A'B}{BD}$,
\item $(A',B',C,D)<(A',B',C',D)$ since
$\frac{C'D}{A'C'}>\frac{CD}{A'C}$,
\item $(A',B',C',D')<(A',B',C',D)$ since
$\frac{C'D'}{B'D'}>\frac{C'D}{B'D}$.
\end{enumerate}

The equality is possible only if all points $A',B',C',D'$ lies on
the corresponding lines, which is impossible since, for example,
the point $B'$ lies in $\pi(S_2)$ which is included in the {\em
interior} of the half-plane $\pi(H_2)$, so $B'\not = B$ and the
inequality in (2) is strict.
\end{proof}

\begin{Lem}\label{pattern}
With conditions as above suppose that four points $A",B",C",D"$
lies on a line $L'$ and
\begin{itemize}
\item $A"\in \partial H_1\setminus \pi(H_4)$
\item $C"\in\partial (\pi(H_3))\setminus \pi(H_4)$
\item $D"\in \partial(\pi(H_4))\setminus H_1$
\end{itemize}
Suppose moreover, that $A"B":B"C":C"D"=A'B':B'C':C'D'$. Then
$B"$
lies in the interior of $\pi(H_2)\cap \pi(H_3)$.
\end{Lem}
\begin{figure}
\centerline{\epsfysize=2in\epsffile{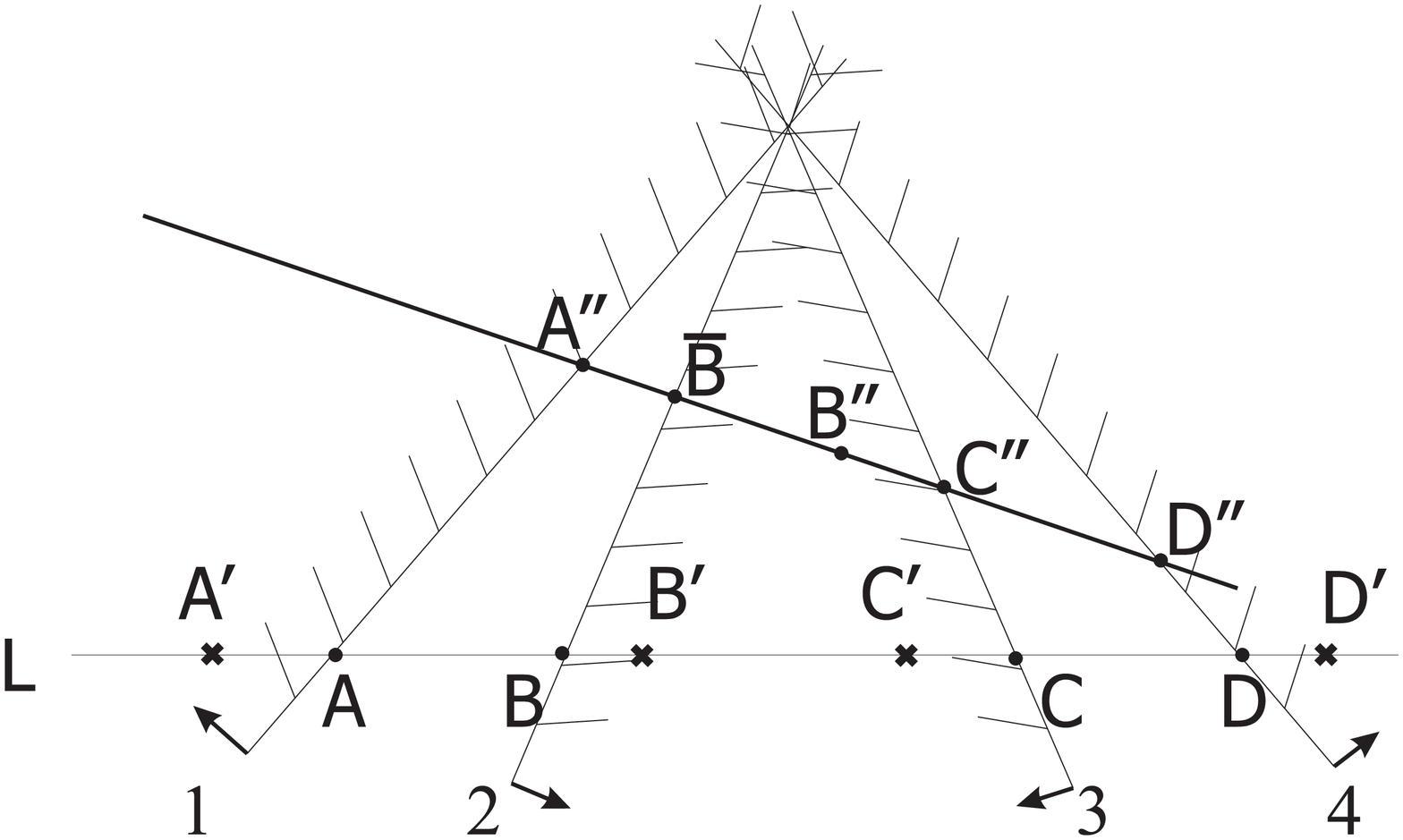}}
\end{figure}

\begin{proof}
This follows directly from the inequality \lemref{double-ratios
relation}. Indeed, let $\overline{B}=L'\cap \partial(\pi(H_2))$.
Then
$(A",\overline{B},C",D")=(A,B,C,D)<(A',B',C',D')=(A",B",C",D")$.
This is equivalent to
$\frac{A"\overline{B}}{\overline{B}D"}<\frac{A"B"}{B"D"}$, which
is possible only if $B"$ is between $ \overline{B}$ and $D"$, i.
e. $B"\in \pi(H_2)$. Since $B"\in [A"C"]$, also $B"\in \pi(H_3)$.
\end{proof}

The Lemma means that the line, which existence is claimed in
\thmref{Restricted Brauer}, can be moved in such a way that it
will still intersect the interior of $H_2$ and  will also
intersect boundaries of $H_1$, $H_3$ and $H_4$.

\subsection{The six non-trivial configurations: contradiction with convex-concavity}

We will call by stencil any  five points $c_1,c_2,c_3,c_4,c_5\in
L_1$ which are projections  points of intersections of some line
$l\subset{\mathbb R}P^3$ with $L_i$, $c_i=l\cap L_i$. Note that
$|c_1c_2|:|c_2c_3|:|c_3c_4|:|c_4c_5|$ is the same for all stencils
and is equal to $h_1:h_2:h_3:h_4$ where $h_i$ are the distances
between $L_i$ and $L_{i+1}$. Evidently this is a necessary and
sufficient condition for five points in $L_1$ lying on a line in
this order to be a projection of points of intersection of some
line in $\mathbb{R}P^3$ with the planes $L_i$.

A projection of a good deformation is a stencil with an additional
property $c_i\in\pi(H_i)$, with at least one of $c_i$ lying in the
interior of $\pi(H_i)$. Vice versa, any such stencil is a
projection of a good deformation.

We can reformulate the \lemref{pattern} using these notations.
\begin{Lem}\label{pattern2}
Suppose that
\begin{enumerate}
\item $H_1\cap \pi(H_4)^c\subset \pi(H_2)^c$,
\item $\pi(H_3)\cap \pi(H_2)\subset \pi(H_4)^c$ and
\item $S_1\cap \pi(H_4)^c\not=\emptyset$.
\end{enumerate}
Then there exists a stencil such that
\begin{enumerate}
\item $c_1\in\partial H_1\cap \pi(H_4)^c$,
\item $c_2$ lies  in the interior of $\pi(H_2)\cap \pi(H_3)$,
\item $c_3\in\partial\pi(H_3)\cap \pi(H_4)^c$ and
\item $c_4\in\partial\pi(H_4)\cap \pi(H_1)^c$.
\end{enumerate}
\end{Lem}

Similar statements hold for all strictly increasing subsequence of
$12345$ consisting of four numbers (i. e. 1245 or 1345 etc.
instead 1234).

\subsubsection{Chebyshev property}

Here we prove that one of  consequences of the Chebyshev property
formulated in \lemref{first in fourth} is that the  set $S_1\cap
\pi(H_4)^c$ in \lemref{pattern} is never empty.

\begin{Lem}\label{trivial-not trivial}
If $S_1\subset \pi(H_4)$ or $S_1\subset \pi(H_5)$  then the
configuration is trivial.
\end{Lem}
\begin{proof}
Indeed, in the first  case $\pi(S_2)$ and $\pi(S_3)$ also lie in
$\pi(H_4)$ by convex-concavity. Indeed, any point of $S_2$ lies on
a segment with endpoints on $S_1$ and $S_4$, and projection of
such a segment  lies entirely in $\pi(H_4)$. The same is true for
$S_3$, so by \lemref{first in fourth} the configuration is
trivial. In the second case $S_i\subset \pi(H_5)$ for
$i=1,2,3,4,5$ and again by \lemref{first in fourth} the
configuration is trivial.
\end{proof}

 In cases C1, C3, C4 and D5 the \lemref{pattern2} and
the \lemref{trivial-not trivial} give immediately existence of a
stencil which is a projection of a good deformation.

\subsubsection{The C1 case}
This is the case $3+2-1+4+5+$. We will consider an equivalent
(after $\b_1^{-2}\b_2$) variant $1+2-3+5-4-$.

If $S_1\cap \pi(H_4)=\emptyset$ then the configuration is trivial
by \lemref{trivial-not trivial}. So $S_1\cap
\pi(H_4)\not=\emptyset$  and the \lemref{pattern2} is
applicable to the subsequence $1+2-3+4-$ of the code.

In the resulting stencil $c_4\in \pi(H_5)$ and $c_1\not\in
\pi(H_5)$. Indeed, the sector $H_1\cap\pi(H_4)^c$ is the smallest
sector bounded by boundaries of half-planes and containing the
point $N$.  Since $N\not\in\pi(H_5)$, so
$H_1\cap\pi(H_4)^c\cap\pi(H_5)=\emptyset$. This means that $-c_4,
c_1\not\in \pi(H_5)$.

Therefore the point $c_5$ of the stencil lies in $\pi(H_5)$.
Therefore the line projecting to this stencil is a good
deformation.

\begin{figure}[h]
\centerline{\epsfysize=1.7in \epsffile{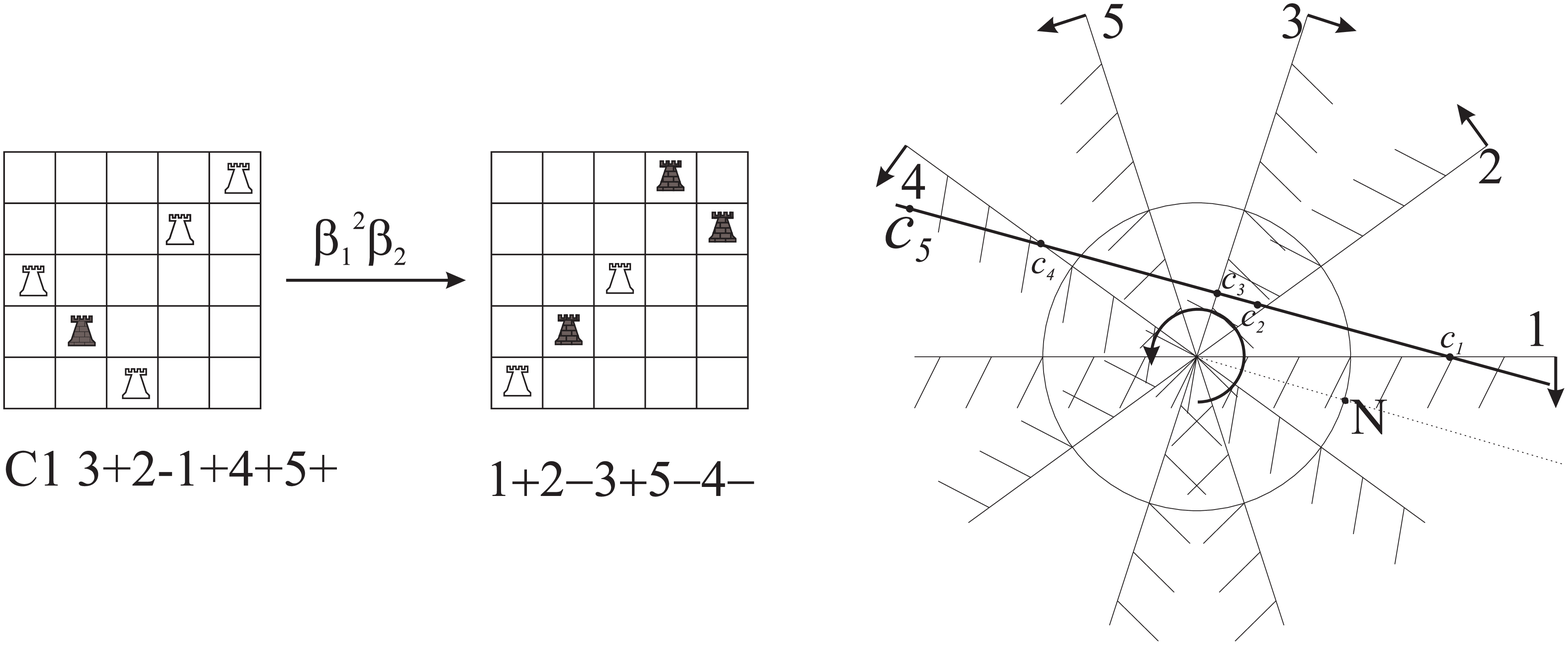}} \caption{The
case C1.}
\end{figure}

\subsubsection{The C3 case}
This is the case of $3+2-4+1+5+$.
We will consider the equivalent (after $\b_2\a_1^6$) case of
$1+2-5-3+4-$.

As above, $S_1\not\subset \pi(H_4)$ by \lemref{trivial-not
trivial}, so  we can apply the  \lemref{pattern2} to the the
subsequence $1+2-3+4-$ of the code, exactly as in the case C1. As
before, $c_1$ lies on $\partial H_1\cap \pi(H_4)^c$ and therefore
in $\pi(H_5)^c$. Also, $c_4\in\partial \pi(H_4)\cap H_1^c$ and
therefore $c_4\in \pi(H_5)$. So $c_5$ also lies in $\pi(H_5)$
since $c_5$ and $c_1$ lie from different sides of $c_4$. Therefore
the stencil given by \lemref{pattern2} is a projection of a good
deformation.
\begin{figure}[h]
\centerline{\epsfysize=1.7in \epsffile{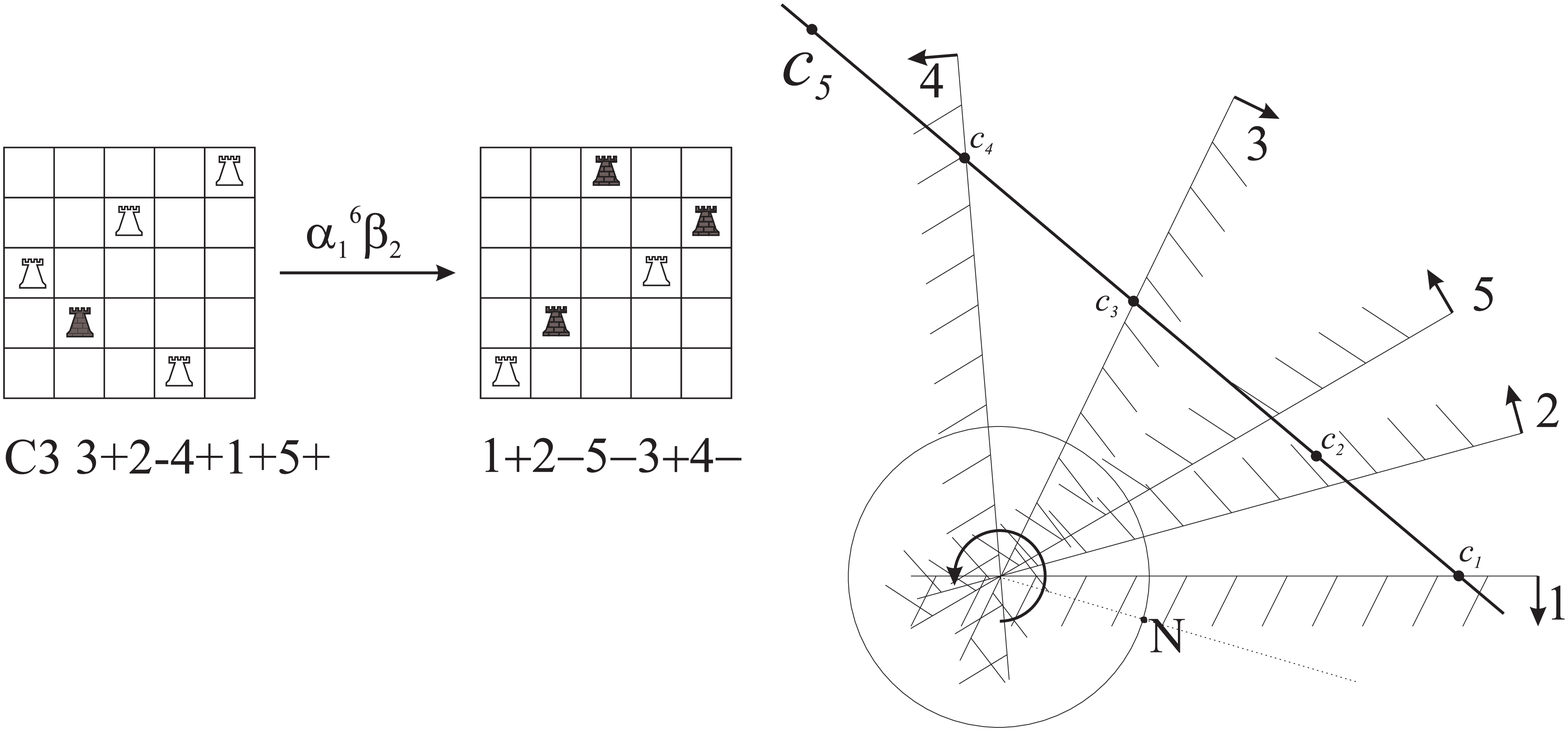}} \caption{The
case C3.}
\end{figure}
\pagebreak

\subsubsection{The C4 case}
This is the case of $3+2-5+1+4+$. We will consider the equivalent
(after $b_1^2 b_2 a_1^{-1}$) case of $1+2-3+5-4+$.

As before, by \lemref{trivial-not trivial}, $S_1\not\subset B_5$.
We apply \lemref{pattern2} to the subsequence $1+2-3+5-$  and get
a stencil with $c_1\in\partial H_1\cap \pi(H_5)^c$,
$c_2\in\pi(H_2)$, $c_3\in\partial \pi(H_3)\cap \pi(H_5)^c$ and
$c_5\in\partial \pi(H_5)\cap H_1^c$. Then $c_4\in\pi(H_4)$.
Indeed, $c_1,c_5\in \pi(H_4)$ and $c_4$ lies between $c_1$ and
$c_5$. So this stencil is a projection of a good deformation.

\begin{figure}[h]
\centerline{\epsfysize=1.7in \epsffile{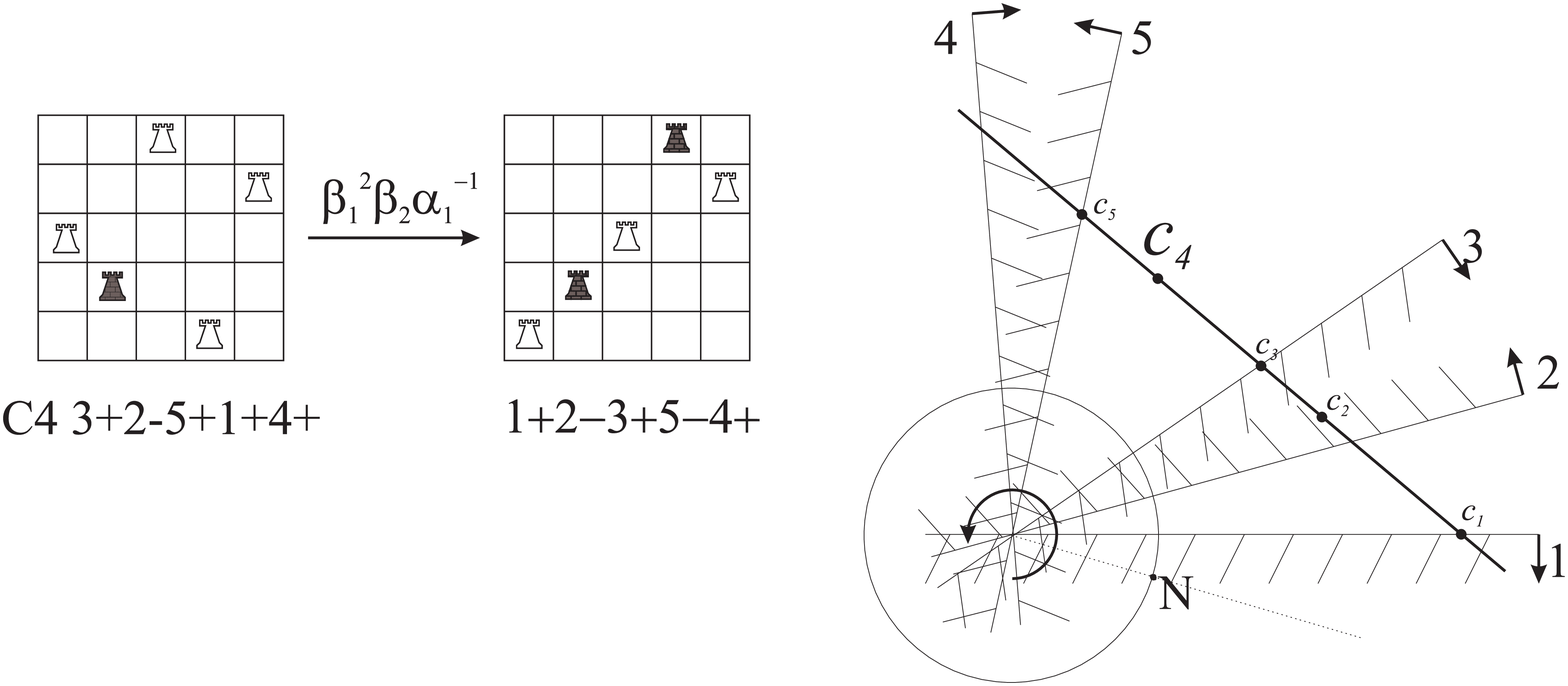}}
\end{figure}

\subsubsection{The D5 case}

This is the case of $4+2-1+3+5+$. It is equivalent (after
$a_1^6b_2$) to the case $1+4-2-3+5-$.

Again,  $S_1\not\subset B_5$ by \lemref{trivial-not trivial}. We
apply \lemref{pattern2} to the subsequence $1+2-3+5-$ and get a
stencil with $c_1\in\partial H_1\cap \pi(H_5)^c$,
$c_2\in\pi(H_2)$, $c_3\in\partial \pi(H_3)\cap \pi(H_5)^c$ and
$c_5\in\partial \pi(H_5)\cap H_1^c$. Now $c_4\in\pi(H_4)$ follows
from the fact that  $c_3,c_5\in \pi(H_4)$ and $c_4$ lies between
$c_3$ and $c_5$. So this stencil is a projection of a good
deformation.
\begin{figure}[h]
\centerline{\epsfysize=1.6in \epsffile{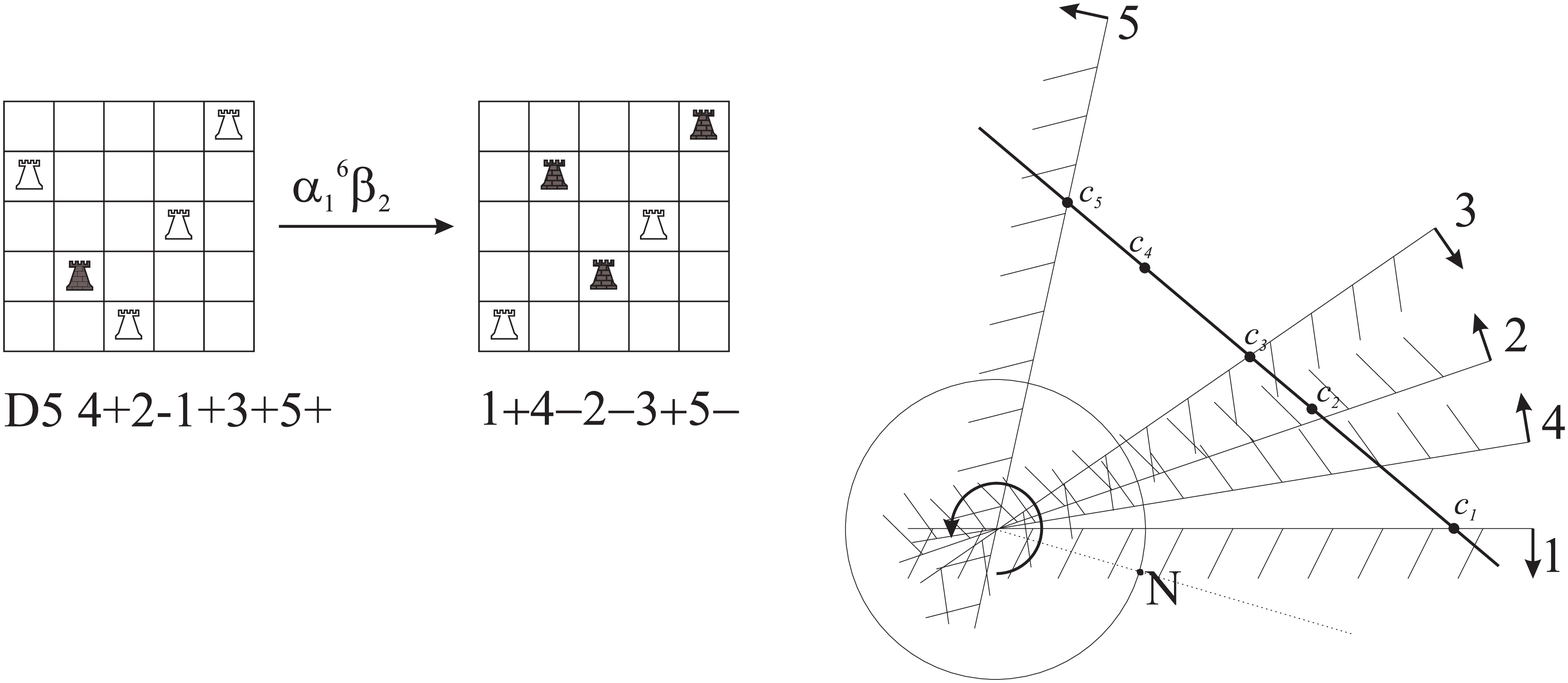}}
\end{figure}
\pagebreak

 In two last cases we should exhibit a little more
inventiveness.

The case C2 requires double application of the
\lemref{pattern2},
whereas in E6 the combinatorial properties of the intersections
contradict to the \thmref{Restricted Brauer}.

\subsubsection{The C2 case}
This is the case of $3+2-1+5+4+$.
After applying $\a_1^3\b_1\b_2$ it will transform to an
equivalent
variant $1+2-3-4+5-$.

By \lemref{trivial-not trivial} $S_1\cap
\pi(H_5)^c\not=\emptyset$. Applying \lemref{pattern2} to the
sequences $1+2-4+5-$ and $1+3-4+5-$ we see that there are two
stencils, one with points $c_1c_2c_3c_4c_5$ and another with
points $c'_1c'_2c'_3c'_4c'_5$, such that the following conditions
hold
\begin{enumerate}
\item $c_1, c'_1\in\partial H_1\cap \pi(H_5)^c$,
\item $c_2\in \pi(H_2)\cap \pi(H_4)$,
\item $c'_3\in \pi(H_3)\cap \pi(H_4)$,
\item $c_4, c'_4\in\partial\pi(H_4)\cap \pi(H_5)^c$ and
\item $c_5, c'_5\in\partial\pi(H_5)\cap \pi(H_1)^c$.
\end{enumerate}

But any two stencils satisfying (1), (4) and (5) differ only by a
dilatation   centered at the origin and these dilatations preserve
$\pi(H_i)$. So $c_3\in \pi(H_3)\cap \pi(H_4)$ and we get the
stencil which is a projection of a good deformation.
\begin{figure}[h]
\centerline{\epsfysize=1.8in \epsffile{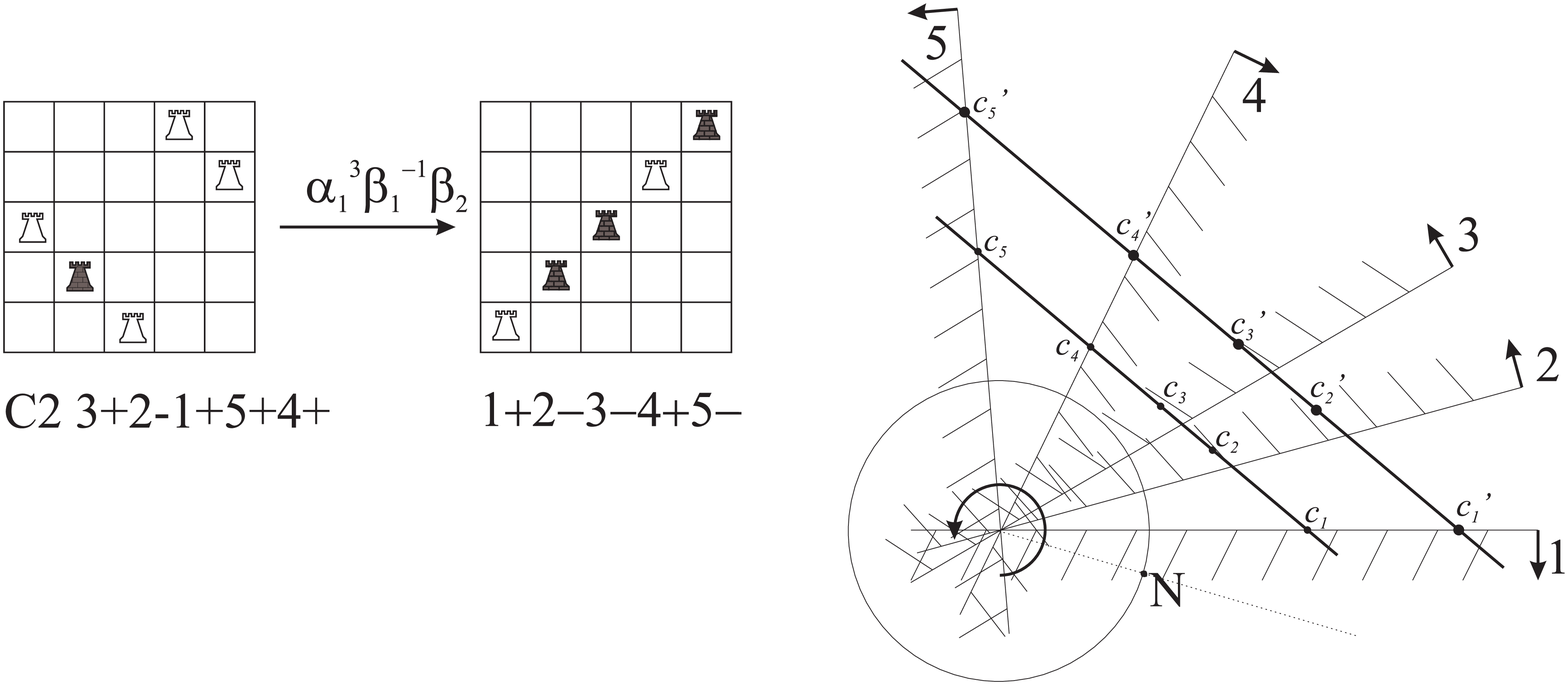}}
\end{figure}

\subsubsection{The E6 case}
This is the case of $4+1+3-5+2+$. It is equivalent (by
$\beta_1^4$) to $1-3+5-2-4+$. Recall that $B_i=\pi^{-1}(H_i)$.

Suppose first that $S_1\cap \pi(H_3)\not=\emptyset$. Similar to
the proof of the \thmref{Restricted Brauer}, we will apply the
\thmref{Browder} to $S_1\cap \pi(H_3)$ as $B$, $S_4\cap B_2$ as
$C$ and $S_2$ and $S_3$ as $A$ and $D$ correspondingly and will
arrive to contradiction.

Construct two mappings, $h_1:CSet(S_1\cap B_3) \to CSet(S_4\cap
B_2)$ and $h_2:CSet(S_4\cap B_2)\to CSet(S_1\cap B_3)$, as in
Namely, take a point $b\in S_1\cap B_3$. There is a line passing
through this point and section $S_2$ and intersecting the section
$S_4$ at point $c$. Since
$\pi(b)\in\pi(S_1)\cap\pi(H_3)\subset\pi(H_2)^c$ and evidently
$\pi(S_2)\subset\pi(H_2)$, we conclude that $\pi(c)\in
\pi(S_4)\cap\pi(H_2)$, i.e. $c\in S_4\cap\pi(H_2)$. The mapping
$h_1$ is the extension to the closed subsets of  $S_1\cap B_3$ of
the mapping sending the points $a$ to the set of all such $c$.
Similarly, to define $h_2$ take any point $c\in S_4\cap B_2$.
There is a line passing through this point and intersecting the
section $S_3$ and the section $S_1$ at a point $a$. Since
$\pi(c)\subset\pi(H_4)\cap\pi(H_2)\subset \pi(H_3)^c$ and
$\pi(S_3)\subset\pi(H_3)$, we get that $a\in S_1\cap B_3$.

In virtue of the \thmref{Browder} this proves existence of a line
intersecting $S_1\cap \pi(H_3)$, $S_2$, $S_3$ and $S_4\cap B_2$.

But this line cannot exist. Indeed, denoting the projections of
the intersection points by $c_1c_2c_3c_4$ we see that  $c_2,c_4\in
\pi(H_4)\cap\pi(H_2)\subset \pi(H_3)^c$ and therefore the point
$c_3$ -- lying between $c_2$ and $c_4$  -- should also belong to
$B_3^c$, which contradicts to $c_3\in \pi(H_3)$.
\begin{figure}[h]
\centerline{\epsfysize=1.8in \epsffile{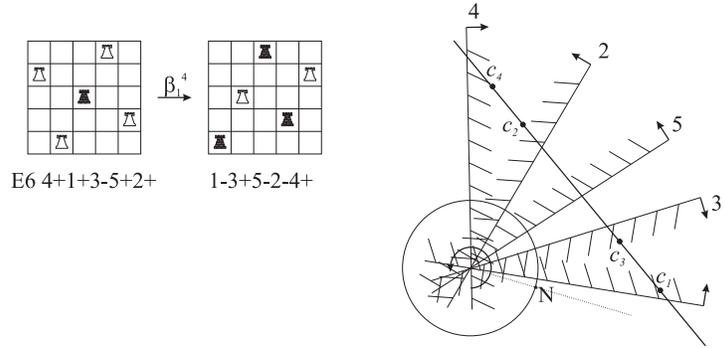}} \caption{The
case of $S_1\cap \pi(H_3)\not=\emptyset$ is impossible.}
\end{figure}

Therefore $S_1\subset \pi(H_3)^c$. By convex-concavity we get that
$\pi(S_4), \pi(H_5) \subset \pi(H_3)$ (any point of these sections
is an endpoint of a segment intersecting $S_3$ with another
endpoint in $S_1$). Therefore $\pi(S_5)\subset \pi(H_5)\cap
\pi(H_3)\subset \pi(H_2)$ and $\pi(S_4)\subset \pi(H_4)\cap
\pi(H_3)\subset \pi(H_2)^c$.

This is incompatible with the existence of lines joining $S_5,
S_4$ and $S_2$ given by convex-concavity condition. Indeed, take
any segment intersecting $S_2,S_4$ and $S_5$ at points $s_2$,
$s_4$ and $s_5$ correspondingly. Its projection $[\pi(s_2),
\pi(s_5)]$ has both ends in $\pi(H_2)$, so $\pi(s_4)\in\pi(H_2)$
as well, which contradicts to
$\pi(s_4)\in\pi(S_4)\subset\pi(H_2)^c$.

\begin{figure}[h]
\centerline{\epsfysize=1.4in \epsffile{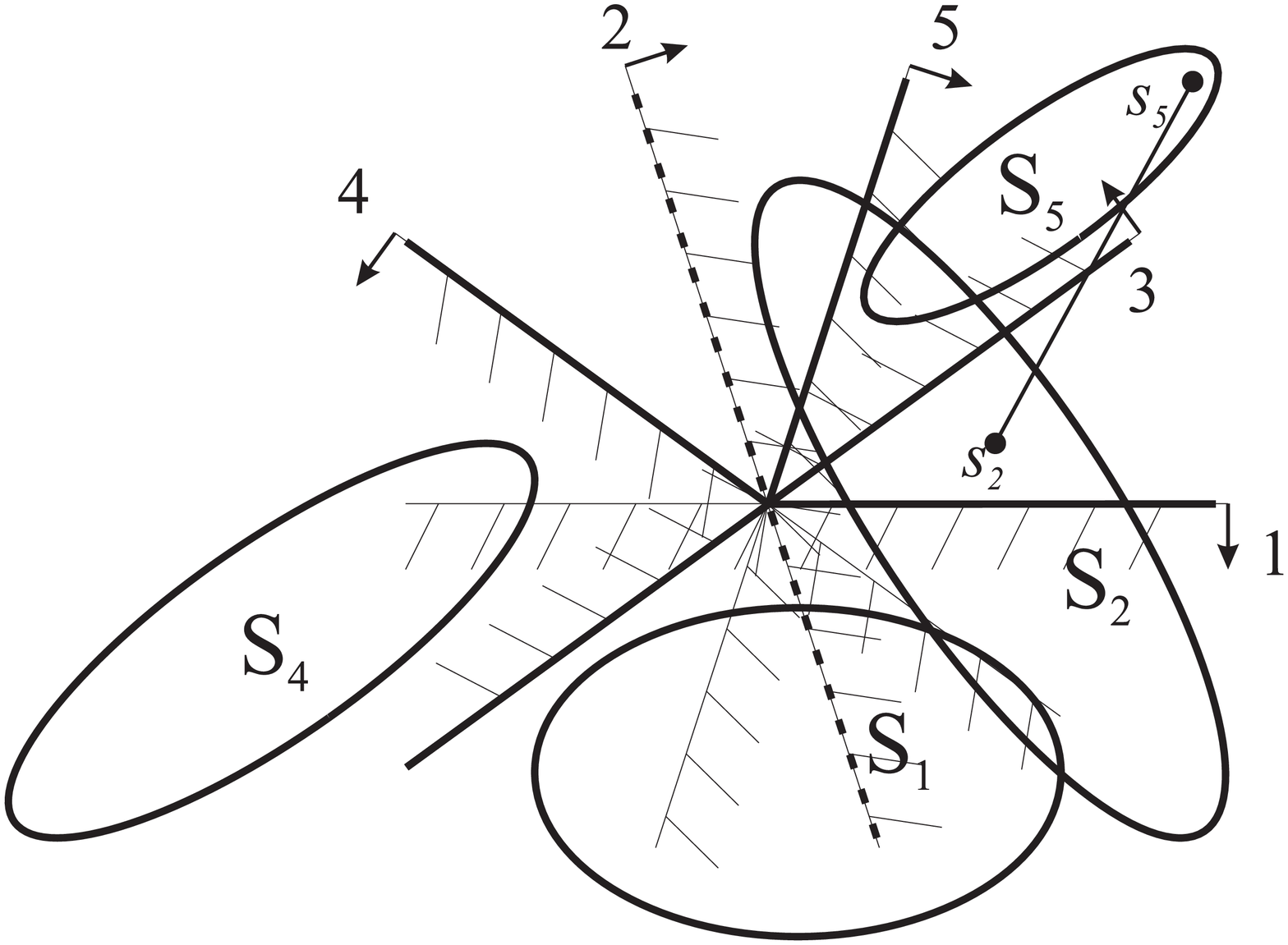}} \caption{The
case of $S_1\subset B_3^c$ is impossible (the previous picture is
rotated by $180^{\circ}$).}
\end{figure}

\bibliographystyle{amsplain}

\end{document}